\documentclass[12pt,a4paper]{amsart}
\usepackage{amsmath}
\usepackage{amscd}
\usepackage{amssymb}

\numberwithin{equation}{section}

\theoremstyle{plain} %% This is the default
\newtheorem{thm}{Theorem}[section]
\newtheorem{cor}[thm]{Corollary}
\newtheorem{lem}[thm]{Lemma}
\newtheorem{prop}[thm]{Proposition}

\newtheorem{defn}[thm]{Definition}

\newtheorem{rem}[thm]{Remark}

\theoremstyle{remark}

\newcommand{\thmref}[1]{Theorem~\ref{#1}}
\newcommand{\propref}[1]{Proposition~\ref{#1}}

\newcommand{\lemref}[1]{Lemma~\ref{#1}}
\newcommand{\corref}[1]{Corollary~\ref{#1}}

\newsavebox{\SmallMathBox}

\def\pdo{\psi{\rm do}}

\def\Ci{C^\infty}

\def\dd{\partial}

\def\Di{D\kern -.65em /}
\def\Dii{D\kern -.45em /}
\def\di{{\dd}\kern -.55em /}
\def\dii{{\dd}\kern -.40em /}

\def\noi{\noindent}

\def\ol{\overline}

\def\to{\rightarrow}
\def\too{\longrightarrow}

\def\re{{\rm Re}}

\def\Bb{{\mathcal B}}
\def\Cc{{\mathcal C}}
\def\Dd{{\mathcal D}}

\def\Nn{{\mathcal N}}

\def\Rr{{\mathcal R}}

\def\={\cong}
\def\>{\supset}
\def\<{\subset}
\def\ii{^{-1}}
\def\si{^{-s}}

\def\12{\frac{1}{2}}
\def\2{\Dd}
\def\3{\Nn}
\def\4{\Rr}
\def\6{\cup}
\def\8{\otimes}
\def\0{^{\circ}}
\def\){\hfill{\ \qed}\enddemo}

\def\a{\alpha}

\def\b{\beta}
\def\Cn{\mathbb{C}^{n}}

\def\d{\delta}

\def\e{\varepsilon}

\def\G{\Gamma}
\def\h{\eta}
\def\th{\theta}

\def\k{\kappa}

\def\la{\lambda}

\def\o{\infty}
\def\s{\sigma}
\def\Si{\Sigma}
\def\z{\zeta}

\def\End{\mbox{\rm End}}

\def\Gr2n{\mbox{${\rm Gr}(\Cn\oplus\Cn)$}}

\def\Grk2n{\mbox{${\rm Gr}_{k}(\Cn\oplus\Cn)$}}
\def\Grk{\mbox{${\rm Gr}_{k}$}}
\def\Gr{\mbox{${\rm Gr}$}}

\def\h{\mbox{\rm h}}

\def\res{\mbox{\rm res}}

\def\Ker{\mbox{\rm Ker}}

\def\Si{S\kern -.65em /}

\def\supp{\mbox{\rm supp}}

\def\tr{\mbox{\rm tr\,}}
\def\Tr{\mbox{\rm Tr\,}}

\def\Vol{\mbox{\rm vol}}

\def\Cf{\mathbb{C}}

\def\Nf{\mathbb{N}}

\def\Rf{\mathbb{R}}

\def\Zf{\mathbb{Z}}

\def\ord{\mbox{\rm ord}}

\def\as{\textsf{a}}
\def\bs{\textsf{b}}
\def\ps{\textsf{p}}
\def\qs{\textsf{q}}
\def\rs{\textsf{r}}

\def\Ds{\textsf{D}}

\def\Is{\textsf{I}}

\def\Ks{\textsf{K}}

\def\Qs{\textsf{Q}}

\def\Ss{\textsf{S}}

\begin{document}

\title[]{{\large The Residue Determinant}\vskip 7mm {\Small {\rm Simon
Scott}}}

\maketitle

%\thispagestyle{empty}

%*****************************************************************

\vskip 5mm

\begin{center}
\noi {\bf Abstract}
\end{center}

\vskip 2mm  The purpose of this paper is to present the
construction of a canonical determinant functional on elliptic
pseudodifferential operators ($\pdo$s) associated to the
Guillemin-Wodzicki residue trace.

The resulting residue determinant functional is multiplicative, a
local invariant, and not defined by a regularization procedure.
The residue determinant is consequently a quite different object
to the zeta function determinant, which is non-local and
non-multiplicative. Indeed, the residue determinant does not arise
as the derivative of a trace on the complex power operators, and
does not depend on a choice of spectral cut. The identification of
a certain residue determinant with the index of an elliptic $\pdo$
shows the residue determinant to be topologically significant.
\footnote{This work arose following conversations with Steve
Rosenberg concerning higher Chern-Weil invariants, my thanks to
him for his support and interest. I am also indebted to Kate
Okikiolu for a helpful suggestion, the essential role of
\cite{Ok1,Ok2} in the current work is evident. I am grateful to
Gerd Grubb for helpful comments and for pointing out a number of
technical improvements. My thanks to Sylvie Paycha for interesting
discussions and, in particular, for pointing out the fact in
Remark(1.11).}

\newpage

\section{Definition and Properties of the Residue Determinant}
Let $A$ be a $\pdo$ of order $\a\in\Rf$ acting on the space of
smooth sections $\Ci(E)$ of a rank $N$ vector bundle $E$ over a
compact boundaryless manifold $M$ of dimension $n$. This means
that in each local trivialization of $E$ with $U\times \Rf^N$,
with $U$ an open subset of $M$ identified with an open set in
$\Rf^n$, and smooth functions $\phi,\psi$ with ${\rm supp}(\phi),
\ {\rm supp}(\psi) \subset U$, then for $x\in U$ and
$f\in\Ci_{c}(U, \Rf^N)$
\begin{equation}\label{e:Alocal}
   (\phi A\psi) f(x) = \frac{1}{(2\pi)^n}\int_{\Rf^n}\int_U e^{i(x-y).\xi}
    \ \as(x,\xi)f(y) \ dy\,d\xi \ ,
\end{equation}
where $\as\in S^{\a}(U)$. We may write $A = {\rm OP}(\as)$ on $U$.

Here, $S^{\a}(U)$ is the symbol space of functions
$\as(x,\xi)\in\Ci(U\times\Rf^n, (\Rf^N)^{*}\otimes\Rf^N)$ with
values in $N\times N$ matrices such that for all multi-indices
$\mu, \nu \in \mathbb{N}^n$,
$\dd_x^{\mu}\dd_{\xi}^{\nu}\as(x,\xi)$ is $O( \ (1 + |\xi|)^{\a -
|\nu|} \ )$, uniformly in $\xi$, and, on compact subsets of $U$,
uniformly in $x$. Write $S(U)$ for $ \cup_{\a\in\Rf}S^{\a}(U)$ and
$S^{-\o}(U)$ for $\cap_{\a\in\Rf}S^{\a}(U)$. Symbols $\as,\bs\in
S(U)$ are said to be equivalent if $\as - \bs \in S^{-\o}(U)$,
written $\as\sim\bs$.

A symbol $\as\in S^{\a}(U)$ is classical (1-step polyhomogeneous)
of degree $\a$ if there is a sequence $\as_0, \as_1, \as_2,
\ldots$ with $\as_j\in\Ci(U\times\Rf^n \backslash \{0\},
(\Rf^N)^{*}\otimes\Rf^N)$ homogeneous in $\xi$ of degree $\a - j$
for $|\xi| \geq 1$ such that $\as(x,\xi) \sim \sum_{j=0}^{\o}
\as_j(x,\xi);$ thus, $\as_j(x,t\xi) = t^{\a - j}\as_j(x,\xi)$ for
$t \geq 1, |\xi| \geq 1$, and
\begin{equation*}
\as(x,\xi) - \sum_{j=0}^{J-1} \as_j(x,\xi) \in S^{\a-J}(U) \ .
\end{equation*}
We may then write $\as\sim(\as_0, \as_1, \ldots)$. A symbol
$\bs\in S(U)$ is called {\em logarithmic} of type $c\in\Rf$ if it
has the form
$$\bs(x,\xi) \sim c\,\log [\xi] \,I + \qs(x,\xi) \ ,$$
where $\qs\sim(\qs_0, \qs_1, \ldots)\in S^0(U)$ is a degree $0$
classical symbol, and $[ \ ] :\Rf^n\to \Rf_+$ is a strictly
positive function with $[\xi]  = |\xi|$ for $|\xi|\geq 1$.

A $\pdo$ $A$ on $\Ci(E)$ is classical of degree $\a$ ({\em resp}.
logarithmic of type $c$) if the local symbol of $A$ in each local
trivialization of $E$ is classical of degree $\a$ ({\em resp}.
logarithmic of type $c$). A logarithmic $\pdo$ has order $\e$ for
any $\e>0$. We denote the space of classical $\pdo$s of order $\a$
({\em resp} less than $\a$) by $\Psi^{\a}(E)$ ({\em resp.
}$\Psi^{<\a}(E)$), and the algebra of all integer order classical
$\pdo$s by $\Psi^{\Zf}(E)$.

The various homogeneous terms $\as_j(x,\xi)$ ({\em resp}
$\qs_j(x,\xi))$ in the local symbol of a classical ({\em resp}.
logarithmic) $\pdo$ do not, in general, have a global invariant
meaning as bundle endomorphisms over $T^*M$. However, it was
observed by Guillemin \cite{Gu} and Wodzicki \cite{Wo2} for
classical $\pdo$s, and extended to logarithmic operators by
Okikiolu \cite{Ok2}, that if $\s(A)_{-n}(x,\xi)$ is the term of
homogeneity $-n$ (so $\s(A)_{-n}(x,\xi) = a_{\a+n}(x,\xi)$ if $A$
is classical of degree $\a$, while if $A$ is logarithmic
$\s(A)_{-n}(x,\xi) = q_{n}(x,\xi)$), then
\begin{equation*}
\frac{1}{(2\pi)^n}\left(\int_{|\xi|=1} \s(A)_{-n}(x,\xi) \ d
S(\xi)\right)\,dx
\end{equation*}
with $d S(\xi)$ the sphere measure on $S^{n-1}$, defines a global
density on $M$. The number
\begin{equation}\label{e:restrace}
  \res(A) = \frac{1}{(2\pi)^n}\int_M\int_{|\xi|=1}
  \tr(\,\s(A)_{-n}(x,\xi)\,)\,
d S(\xi)\, dx
\end{equation}
is the Guillemin-Wodzicki {\em residue trace} of the $\pdo$ $A$.
Evidently, if $A$ is a classical $\pdo$ of order $\a$, then
\begin{equation}\label{e:restracezero}
\a \notin \Zf  \hskip 7mm {\rm or} \hskip 7mm \a < -n \ \ \
\Rightarrow  \ \ \ \res(A) = 0 \ ,
\end{equation}
and $\res$ drops down to a map on the quotient algebra
\begin{equation}\label{e:restracequotient}
  \res : \Psi^{\Zf}(E)/\Psi^{-\o}(E) \too \Cf \ .
\end{equation}
The following linearity properties of the residue trace are
immediate from its definition.

\begin{lem}\label{lem:linearity}
Let $A, B$ be $\pdo$s. If $A$ and $B$ are both logarithmic, or, if
$A$ is classical of order $\a\in\Zf$ and $B$ is logarithmic, or,
if $A$ is classical of order $\a\in\Rf$ and $B$ is classical of
order $\b\in\Rf$ such that\footnote{If $A$ and $B$ are classical
and $\a - \b \notin \Zf$ then $A + B$ is not a classical $\pdo$
(the symbol is then not 1-step, its expansion does not drop in
integer orders). Consequently, though $\Psi^{\Zf}(E)$ is an
algebra, the space $\Psi(E)= \cup_{\a\in\Rf} \Psi^{\a}(E)$ is not,
but forms a semi-group with respect to the usual composition
product. This is relevant for linearity properties of traces on
classical $\pdo$s, and the reason why $\s(\log A)(x,\xi) =
-(d/ds)|_{s=0}\s(A^{-s})(x,\xi)$, as a limit of differences of
classical symbols, is not quite classical.} $\a - \b \in \Zf \ ,$
then
\begin{equation}\label{e:restracelinear}
  \res(A + B) = \res(A) + \res(B) \ .
\end{equation}
\end{lem}

\vskip 2mm

The characterizing tracial property of $\res$ is due to Guillemin,
Wodzicki and extended to include logarithmic operators by
Okikiolu:

\begin{prop}\label{p:restrace}\cite{Wo2},\cite{Gu},\cite{Ok2}.
Let $A, B$ be classical or logarithmic $\pdo$s. Then
\begin{equation}\label{e:restracial}
  \res([A,B]) = 0 \ .
\end{equation}
\end{prop}

\vskip 1mm

\noi It follows that \eqref{e:restracequotient} is a trace
functional. It is, moreover, projectively unique.

\vskip 1mm

This has consequences for determinants.

\vskip 1mm

Let $A$ be a classical $\pdo$ $A$ of order $\a$ admitting a
principal angle $\th$, meaning that the principal symbol
$\as_0(A)(x,\xi)$ considered as a bundle endomorphism over
$T^*X\backslash 0$ has no eigenvalue on the spectral cut $R_{\th}
= \{re^{i\th} \ | \ r\geq 0\}$; in particular, $A$ is elliptic.
Then, as recalled below, the functional calculus constructs the
log symbol $(\log_{\th}\as)_{-n}(x,\xi)$ of homogeneity $-n$ and
\begin{equation*}
\frac{1}{(2\pi)^n}\left(\int_{|\xi|=1} (\log_{\th}\as)_{-n}(x,\xi)
\ d S(\xi)\right)\,dx
\end{equation*}
defines a global density on $M$ \cite{Ok2}. If $\a>0$, then
$\log_{\th}A$ exists as a logarithmic $\pdo$ of type $\a$ and
$(\log_{\th}\as)_{-n}(x,\xi) = \s(\log_{\th}A)_{-n}(x,\xi)$. We
can hence define canonically the following determinant functional
on classical $\pdo$s.
\begin{defn}
The residue determinant $\det_{\rm res}A$ of a classical $\pdo$
$A$ with principal angle $\th$ is the complex number
\begin{equation}\label{e:resdet}
  \log{\rm det}_{\rm res} A  := \res(\log A) \ ,
\end{equation}
that is,
\begin{equation}\label{e:resdet2}
\log{\rm det}_{{\rm res}}A  :=
\frac{1}{(2\pi)^n}\int_M\int_{|\xi|=1}\tr\left((\log_{\th}\as)_{-n}(x,\xi)\right)
\ d S(\xi)\, dx \ .
\end{equation}

\vskip 2mm

\end{defn}
\begin{rem}
From \lemref{lem:linearity}, one has the linearity
\begin{equation}\label{e:restracelinear2}
\res(\log A + \log B) = \res(\log A) + \res(\log B) \
\end{equation}
for all classical $\pdo$s $A$ and $B$, of any real orders
$\a,\b\in\Rf$.
\end{rem}

\vskip 3mm

The properties of the residue determinant are as follows.

\vskip 3mm

\begin{prop}\label{p:localnotheta}
The residue determinant ${\rm det}_{\res} A $ is a local
invariant, depending only on the first $n+1$ homogeneous terms
$\as_0, \as_1, \ldots, \as_n$ in the local symbol of $A$, and is
independent of the choice of principal angle $\th$ used to define
$\log_{\th}A$.
\end{prop}

\vskip 3mm

\noi The subscript $\th$ is therefore omitted in the notation
\eqref{e:resdet}.

\vskip 1mm

\begin{cor}\label{cor:detresstable}
Let $A\in \Psi^{\a}(E), \ S\in \Psi^{\s}(E)$ with $\a  > \s + n$.
Then
\begin{equation}\label{e:detresstable}
{\rm det}_{{\rm res}}(A + S) = {\rm det}_{{\rm res}}(A) \ .
\end{equation}

\end{cor}

\vskip 3mm

The residue determinant is multiplicative:

\vskip 3mm

\begin{thm}\label{t:mult} Let $A, B$ be classical
$\pdo$s of order $\a ,\b\in\Rf$, and suppose that $A, B, AB$ admit
principal angles. Then
\begin{equation}\label{e:detresmult} {\rm
det}_{{\rm res}}(AB) = {\rm det}_{{\rm res}}(A)\cdot {\rm
det}_{{\rm res}}(B) \ .
\end{equation}
\end{thm}

\vskip 2mm

\noi The residue determinant does not vanish on non-invertible
operators, however.

\vskip 1mm

In fact, since an elliptic operator is invertible modulo smoothing
operators, the above properties imply that it never vanishes.

\vskip 3mm

If $A\in\Psi(E)$ has order $\a > 0$ and is invertible, then
$\z(A,0)|^{{\rm mer}}$, the meromorphically continued spectral
zeta-function of $A$ evaluated at $s=0$ (see Sect.2), is known to
have the properties in \propref{p:localnotheta}; the relation with
$\det_{\rm res}A$ is the following.

\vskip 3mm

\begin{thm}\label{t:resdetzetazero}
Let $A$ be a classical $\pdo$ of order $\a>0$ with principal
angle. If $A$ is invertible
\begin{equation}\label{e:resdetzetazero}
{\rm det}_{{\rm res}}(A) = e^{-\a \,\z(A,0)|^{{\rm mer}}} \ .
\end{equation}
\vskip 1mm \noi When $A$ is not invertible
\begin{equation}\label{e:resdetzetazero1}
\hskip 12mm {\rm det}_{{\rm res}}(A) = e^{-\a \,(\z(A,0)|^{{\rm
mer}} \, +  \ {\rm h}_0(A))} \ ,
\end{equation}\vskip 1mm
\noi where $\h_0(A) = \Tr(\Pi_0(A))$ and $\Pi_0(A)$ is a
projection  onto the finite-dimensional generalized $0$-eigenspace
$E_0(A) = \{\tau\in\Ci(E) \ | \ A^N\tau = 0 \ {\rm some} \ N\in
\Nf\}$ (the projection $\Pi_0(A)$ is defined in
\eqref{e:powers2}). If $\Ker(A^2) = \Ker(A)$, in particular for A
self-adjoint, then $\h_0(A) =\dim \Ker(A)$.
\end{thm}

\vskip 2mm

The additional ${\rm h}_0(A)$ term on the right-side of
\eqref{e:resdetzetazero1} thus corrects for the discontinuities of
$\z(A,0)|^{{\rm mer}}$ at non-invertible $A$.  Notice, also, that
$\z(A,0)|^{{\rm mer}} $ is locally determined only if $\h_0(A) =
0$, otherwise it is $\z(A,0)|^{{\rm mer}} \, + \ {\rm h}_0(A)$
which is local; this follows from \propref{p:localnotheta} and
\eqref{e:resdetzetazero}, and is seen directly in the proof.

\vskip 2mm

\begin{rem}
\noi {\rm [1]} \ The number $\z(A,0)|^{{\rm mer}} :=
\Tr(I.A\si)|^{{\rm mer}}_{s=0}$ defines a quasi- (or weighted-)
trace of the identity operator $I$ on $\Ci(E)$ and hence
\eqref{e:resdetzetazero} associates $\det_{\rm res}$ with a notion
of regularized dimension, rather than regularized volume.

\vskip 2mm

\noi {\rm [2]} \ The continuity of $\det_{\rm res}$ on families of
admissible operators in $\Psi(E)$ contrasts with the
$\zeta$-determinant, which is continuous only on families of
invertible operators.

\end{rem}

\begin{rem}
\noi Subsequently two proofs of \eqref{e:resdetzetazero1} of
independent interest have been given. In \cite{Pa} (see also
\cite{PaSc}) the identity is proved using a microlocal result of
\cite{KV}. In \cite{Gr2} a proof is obtained via the resolvent
trace $\Tr((A-\la I)^{-k})$.
\end{rem}

\begin{rem}Since $\log_{\th}A$ is `almost' in the subalgebra
$\Psi^{\leq 0}(E)\cap \Psi^{\Zf}(E)$ on which the residue trace is
not the unique trace \cite{PaRo}, ${\rm det}_{{\rm res}}$ is not
quite the unique multiplicative functional  on elliptic $\pdo$s.
Indeed ${\rm det}_0$ defined by  $\log {\rm det}_0 (A) = ({\rm
vol}_{\sigma}(S^{*}M))\ii\int_{S^{*}M} \log {\rm det}(a_0(x,\xi))
\ d \s$, where $a_0(x,\xi)$ is the leading symbol of $A$ and $d\s$
is a volume form on the cosphere bundle $S^* M$, is
multiplicative.
\end{rem}

\vskip 5mm

\noi {\bf Example.} \ \ Let $\Sigma$ be a closed Riemann surface
and $E$ a complex vector bundle of degree ${\rm deg}(E) =
\int_{\Sigma} c_1 (E)$. Let $\ol\partial_{\Sigma} : \Ci(E) \too
\Ci(E\otimes T^{0,1}\Sigma)$ be an invertible
$\ol\partial$-operator; thus locally, $\ol\partial_{\Sigma}
=(\partial_{\overline{z}} + a(z)) d\overline{z}$. Then, from
\eqref{e:resdetzetazero} and \cite{Gi} Thm(4.1.6)(see also
\cite{Bo} \S1.5), we have
\begin{equation}\label{e:Dbaronsurfaces}
{\rm det}_{{\rm res}}(\ol\partial_{\Sigma}^{ \, *}
\ol\partial_{\Sigma}) = \exp\left(- \, {\rm deg}(E) -
\frac{\chi(\Sigma) \, {\rm rk}(E)}{3}\right) \ ,
\end{equation}
with $\chi(\Sigma)$ the Euler number, and ${\rm rk}(E)$ the rank
of $E$.

\vskip 1mm

 This is independent of $\ol\partial_{\Sigma}^{ \, *}
\ol\partial_{\Sigma}$, but in general the residue determinant of a
second-order differential operator over a surface will depend on
the complete symbol. In the case of an invertible operator
$\Delta_g$ of Laplace-type one has for $t\in\Rf$
\begin{equation}\label{e:laplacianonsurfaces}
{\rm det}_{{\rm res}}(\Delta_g + tI) = \exp\left(
\frac{A_g(\Sigma) \, {\rm rk}(E)}{2\pi}\ t -
\frac{1}{2\pi}\int_{\Sigma}\tr(\e_x(\Delta_g))\, dx -
\frac{\chi(\Sigma) \, {\rm rk}(E)}{3}\right) \ ,
\end{equation}

\noi where $A_g(\Sigma)$ is the surface area of $\Sigma$ with
respect to a Riemannian metric $g$, and $\e(\Delta_g)$ is the
unique element of $\Ci(\End(E))$ and $\nabla$ the unique
connection such that $ \Delta_g =
-\sum_{i,j}g^{ij}(x)\nabla_i\nabla_j + \e_x(\Delta_g)$. On the
other hand, if $(M,g)$ is a 4-manifold and $\Delta_g$ the
Laplace-Beltrami operator then
\begin{equation}\label{e:laplacianonon4manifolds}
{\rm det}_{{\rm res}}(\Delta_g + tI) = \exp\left( -\frac{
\Vol_g(M)}{16\pi^2}\ t^2 + \frac{1}{24\pi^2}\int_M \k_M \, dx \ t
\right) \, {\rm det}_{{\rm res}}(\Delta_g) \ ,
\end{equation}
with $\Vol_g(M)$ the Riemannian volume, $\k_M$ the scalar
curvature. More generally, if $(M,g)$ is a Riemannian manifold of
dimension $2m$, $E$ a vector bundle, and $\Delta_g$ an operator on
$\Ci(E)$ of Laplace-type, meaning $\Delta_g$ is a second-order
differential operator with scalar principle symbol
\begin{equation}\label{e:laplacetype}
\s(\Delta_g)_{\,2}(x,\xi) = |\xi|_{g(x)}^2 \,I :=
-\sum_{i,j=1}^{2m}g^{ij}(x)\xi_i\xi_j \,I \ ,
\end{equation}
where $\xi = (\xi_1,\ldots,\xi_{2m})\in\Rf^{2m}$, then for $\e
>0$ the heat operator $e^{-\e\Delta_g}$ is a smoothing operator
with heat trace expansion as $\e\to 0+$
\begin{equation}\label{e:heattrace}
\Tr(e^{-\e\Delta_g}) = \frac{c_{-m}(\Delta_g)}{\e^m} + \ldots +
\frac{c_{-1}(\Delta_g)}{\e} + c_0(\Delta_g) + O(\e)
\end{equation}
with locally determined coefficients $c_j(\Delta_g)$;
specifically, \eqref{e:laplacetype} easily implies
\begin{equation}\label{e:cm}
c_{-m}(\Delta_g) = \frac{\Vol_g(M)\,{\rm rk}(E)}{(4\pi)^m} \ ,
\end{equation}
while, by standard transition formulae (see for example \cite{GS},
\cite{Gr}), \eqref{e:resdetzetazero1} becomes
\begin{equation}\label{e:resdetlaplace}
{\rm det}_{{\rm res}}(\Delta_g) = e^{- 2 \,c_0(\Delta_g)} \ .
\end{equation}
(If $M$ is odd-dimensional ${\rm det}_{{\rm res}}(\Delta_g) =1$.)
On the other hand, for $t\in\Rf$ one has $$\Tr(e^{-\e(\Delta_g +
tI)}) = e^{-\e \, t}\,\Tr(e^{-\e\Delta_g}) \ ,$$ and so from
\eqref{e:heattrace}
\begin{equation}\label{e:zetazerot}
c_0(\Delta_g + tI) = \frac{(-1)^m}{m!}\, c_{-m}(\Delta_g) \ t^m +
\frac{(-1)^{m-1}}{(m-1)!}\, c_{-m+1}(\Delta_g) \ t^{m-1} + \ldots
+ c_0(\Delta_g) \ .
\end{equation}

\noi Thus \eqref{e:cm} and \eqref{e:resdetlaplace} yield

\begin{eqnarray*}
{\rm det}_{{\rm res}}(\Delta_g + t I) & = &
\exp\left(\frac{2(-1)^{m+1}}{ (4\pi)^{m} \, m!}\Vol_g(M)\,{\rm
rk}(E)
\ t^m  \ \ + \right. \\
& & \hskip 15mm \left.
2\sum_{j=1}^{m-1}\frac{(-1)^{m-j+1}}{(m-j)!} \, c_{-m+j}(\Delta_g)
\ t^{m-j}  \right) \, {\rm det}_{{\rm res}}(\Delta_g) \ ,
\end{eqnarray*}

\noi the specific formulas \eqref{e:laplacianonsurfaces},
\eqref{e:laplacianonon4manifolds} now following from \cite{Gi}
Thm(4.1.6) and \cite{McSi}.

\begin{rem}
Let $\omega\in\Ci(M)$ and let $g_{\omega} = e^{2\omega}g$. A
consequence of \eqref{e:resdetlaplace} is that the Generalized
Polyakov Formula of \cite{BrOr} for the relative zeta-determinant
of $\Delta_{g_{\omega}}$ may be stated as
\begin{equation*}
\log\frac{{\rm det}_{\zeta}(\Delta_{g_{\omega}})}{{\rm
det}_{\zeta}(\Delta_{g_0})} = \int_0^1 \log {\rm det}_{{\rm
res}}(\Delta_{g_{\epsilon\omega}}) \ d\epsilon \ .
\end{equation*}
\end{rem}

\vskip 4mm

\noi Turning matters around, \eqref{e:resdetzetazero1} can be used
to deduce properties of $\z(A,0)|^{{\rm mer}}$.

\vskip 1mm

Since the proof of \thmref{t:resdetzetazero} demonstrates that the
equality in \eqref{e:resdetzetazero1} holds logarithmically
\begin{equation}\label{e:resdetzetazero2}
{\rm res}(\,\log A \,) = -\a \,\left(\z(A,0)|^{{\rm mer}} +
\h_0(A)\right)\ ,
\end{equation}

\vskip 2mm

\noi (in fact, it holds pointwise on $M$ as an equality between
densities as shown in the proof of \thmref{t:resdetzetazero},)
then, combined with \eqref{e:detresmult} which also holds
logarithmically, (but not as a pointwise identity of densities) we
have the following.

\vskip 2mm

\begin{cor} \  Let $A$, $B$ be classical $\pdo$s, admitting
principal angles, and which have positive orders $\a
,\b\in\Rf_{+}$. Then the function $Z(\s(A)) := -\a
\,\left(\z(A,0)|^{{\rm mer}} + \h_0(A)\right)$ is additive:
\begin{equation}\label{e:zetaatzeroadditive}
Z(\s(AB)) = Z(\s(A)) + Z(\s(B))\ .
\end{equation}
\end{cor}

\vskip 3mm

\begin{rem} \noi {\rm [1]} The additivity
\eqref{e:zetaatzeroadditive} is referred
to in \cite{Ka} and in the introduction to \cite{KV} (whose
notation is respected here) as a property known to Wodzicki,
though no proof appeared. On the other hand, Wodzicki defined in
\cite{Wo3} a determinant for order zero $\pdo$s path connected to
the identity which, for such operators, coincides with ${\rm
det}_{{\rm res}}$.

\vskip 3mm

\noi {\rm [2]} In \cite{PaSc} it is shown that
\eqref{e:resdetzetazero2} extends as an exact splitting formula
into local and global components of the generalized zeta function
$\z_{\th}(A,Q,s)|^{{\rm mer}} = \Tr(AQ\si)|^{{\rm mer}}$. On the
other hand, Grubb \cite{Gr2} has recently shown using resolvent
methods that \eqref{e:resdetzetazero2} extends to certain classes
of boundary value problems.
\end{rem}

\vskip 1mm

Conversely, \eqref{e:resdetzetazero1},
\eqref{e:zetaatzeroadditive} combine to prove \eqref{e:detresmult}
when $\a, \b \in \Rf_{+}$.

\vskip 2mm

From the Atiyah-Bott-Seeley $\z$-function index formula (see for
example \cite{Sh}), a further immediate consequence of
\eqref{e:resdetzetazero} is the following `super' residue
determinant formula for the index ${\rm Index}(\Ds) =
\dim\Ker(\Ds) - \dim\Ker(\Ds^*)$ of a general elliptic operator
$\Ds:\Ci(E^+)\too\Ci(E^-)$ of order $d > 0$.

\begin{cor}
\begin{equation}\label{e:resdetindex}
\frac{{\rm det}_{{\rm res}}(\Ds^{*} \Ds + I)}{{\rm det}_{{\rm
res}}(\Ds \Ds^{*} + I)} = e^{-2d \, {\rm Index}(\Ds)} \ .
\end{equation}
Equivalently,
\begin{equation}\label{e:resdetindex2}
{\rm Index}(\Ds) = \frac{1}{2d}\left({\rm res}\log(\Ds \Ds^{*} +
I) - {\rm res}\log(\Ds^{*} \Ds + I)\right)
 \ .
\end{equation}

\end{cor}

\vskip 2mm

\begin{rem}
In contrast, there is not a formula for the index using the
residue trace of a classical (non-logarithmic) $\pdo$. The
identities \eqref{e:resdetzetazero1} and  \eqref{e:resdetindex2}
lead to an alternative elementary proof of the local Atiyah-Singer
index theorem \cite{ScZa}.
\end{rem}

\vskip 2mm

For order zero operators we have:

\vskip 1mm

\begin{thm}\label{t:idpluslowerorder}
If $A$ has order $\a =0$ and has the form $A = I + \Qs$ with $\Qs$
a classical $\pdo$ of  negative integer order $ k < 0$, then
\begin{equation}\label{e:degreezero}
{\rm det}_{{\rm res}}(I + \Qs) =
\prod_{j=1}^{\left[\frac{n}{|k|}\right]}e^{\frac{(-1)^{j+1}}{j}\,
\res(\Qs^j)} \ .
\end{equation}
\end{thm}

\vskip 1mm

\begin{rem}
\noi {\rm [1]}  For order zero operators the zeta-function at zero
in \eqref{e:resdetzetazero1} is thus replaced by a relative zeta
function at zero (see \thmref{t:zetapolynomial}, and comments
around \eqref{e:relzeta}).

\vskip 1mm

\noi {\rm [2]} \  Generalizing the Fredholm determinant (p=1)
there is a well-known notion of $p$-determinant $\det_p(I + Q)$
for $Q$ in the $p^{th}$ Schatten ideal $L_p$, which is not
multiplicative for $p>1$, but for $Q\in L_1$ the following formula
holds \cite{Si} in analogy to \eqref{e:degreezero}
\begin{equation*}
    \frac{{\rm det}_p(I+Q)}{{\rm det}_1(I+Q)} =
    \prod_{j=1}^{p-1}e^{\frac{(-1)^{j}}{j}\Tr(Q^j)} \ .
\end{equation*}
\end{rem}

\vskip 2mm

\thmref{t:idpluslowerorder} evidently implies:

\vskip 1mm

\begin{cor}\label{c:detres=1}
\begin{equation}\label{e:det=1}
{\rm det}_{{\rm res}}(I + \Qs) = 1 \hskip 10mm {\rm if } \hskip
4mm \ord(\Qs) < -n \ , \ \ord(\Qs)\in\Zf \ .
\end{equation}
Hence ${\rm det}_{\res}$ drops down to a multiplicative function
on the `determinant Lie group' $\widetilde{G} = \Psi^{*}(E)/(I +
\Psi^{-\o}(E)),$ of Kontsevich-Vishik \cite{KV}, where
$\Psi^{*}(E)$ is the group of invertible elliptic $\pdo$s.
(\eqref{e:det=1} also follows from \corref{cor:detresstable}.)
\end{cor}

\vskip 1mm

\noi {\bf Example}. \ \ Let $\Delta_g$ be an invertible
generalized Laplacian on a closed Riemannian manifold $(M,g)$ of
dimension $2m$. Thus $\Delta_g$ has principal symbol as in
\eqref{e:laplacetype}, and with that notation, we therefore have
$\s(\Delta_g^{-m})_{-2m}(x,\xi) = |\xi|_{g(x)}^{-2m}I$. Whence,
with $S^{2m-1}$ the Euclidean $(2m-1)$-sphere,
\begin{eqnarray}
\res(\Delta_g^{-m})  & = & \int_M \frac{1}{(2\pi)^{2m}} \left(
\int_{|\xi|_{g(x)}=1} \tr(I)\  dS_{g(x)}(\xi) \right) dx  \nonumber \\
& = &  \frac{1}{(2\pi)^{2m}} \int_M \sqrt{{\rm det}(g(x))} \left(
\int_{S^{2m-1}} \ dS(\xi) \right) dx \ {\rm rk}(E) \nonumber \\
& = &  \frac{\Vol(S^{2m-1})}{(2\pi)^{2m}} \, \Vol_g(M) \, {\rm
rk}(E)
 \label{e:resdeltam} \ .
\end{eqnarray}
Since $\Vol(S^{2m-1}) = 2\pi^m / (m-1)!$ we therefore have from
\eqref{e:degreezero}
\begin{equation}\label{e:lap} {\rm det}_{{\rm res}}(I + \Delta^{-m}) =
\exp\left(\frac{\Vol_g(M)\, {\rm rk}(E)}{2^{2m-1} (m-1)!\, \pi^m}
\right) \ .
\end{equation}

\vskip 7mm

\noi {\bf Example}. \  It is instructive to check how the
multiplicativity property of ${\rm det}_{\res}$ works for this
class of $\pdo$s. As a simple case, consider $\pdo$s  $\Qs_1,
\Qs_2$ for which $3\,\ord(\Qs_i) < - n \leq 2\,\ord(\Qs_i) < 0$ --
for example, operators of order -2 on a 4-manifold. Then according
to \thmref{t:idpluslowerorder}
\begin{equation}\label{e:minus2on5}
\log{\rm det}_{\res}(I + \Qs_i) = \sum_{p=1
}^2\frac{(-1)^{p+1}}{p}\ \res(\Qs_i^p)
 =  \res(\Qs_i) -  \frac{1}{2}\,\res(\Qs_i^2) \ .
 \end{equation}
So
\begin{equation}\label{e:ex1a}
\log{\rm det}_{\res}(I + \Qs_1) + \log{\rm det}_{\res}(I + \Qs_2)
\end{equation}
$$ =  \res(\Qs_1) + \res(\Qs_2) -  \frac{1}{2}\,\res(\Qs_1^2) -
\frac{1}{2}\,\res(\Qs_2^2) \ . $$ On the other hand,
$(I+\Qs_1)(I+\Qs_2) =  I+(\Qs_1 + \Qs_2 + \Qs_1 \Qs_2)$ and
\eqref{e:minus2on5} applies with $\Qs_i$ replaced by $\Qs_1 +
\Qs_2 + \Qs_1 \Qs_2$ so that
\begin{equation}\label{e:ex1b}
\log{\rm det}_{\res}((I+\Qs_1)(I+\Qs_2))
\end{equation}
$$=  \res(\Qs_1 + \Qs_2 + \Qs_1 \Qs_2) -
\frac{1}{2}\,\res\left((\Qs_1 + \Qs_2 + \Qs_1 \Qs_2)^2\right) \
.$$ Since the $\Qs_i$ have integer order
$$\res(\Qs_1 + \Qs_2 + \Qs_1 \Qs_2) = \res(\Qs_1) + \res(\Qs_2) + \res(\Qs_1
\Qs_2)$$ while by the tracial property \eqref{e:restracial} of
$\res$
\begin{eqnarray*}
\res\left((\Qs_1 + \Qs_2 + \Qs_1 \Qs_2)^2\right)& = & \res(\Qs_1^2
+ \Qs_2^2 + \Qs_1 \Qs_2 + \Qs_2 \Qs_1 + \ {\rm terms\ of \ order <
-n }) \\ & = &  \res(\Qs_1^2) + \res(\Qs_2^2) + 2\,\res(\Qs_1
\Qs_2) \ .
\end{eqnarray*}
 Hence
\eqref{e:ex1a} and \eqref{e:ex1b} are equal.

\vskip 4mm

An application of properties in the previous theorems yields the
following formula.

\vskip 3mm

\begin{thm}\label{t:zetapolynomial}
Let $A \in\Psi^{\a}(E), B_0\in\Psi^{\b_0}(E), \ldots , B_d
\in\Psi^{\b_d}(E)$ be classical $\pdo$s. Assume  $\a = \ord(A)> 0$
and $ \a - \b_j \in \Nf\backslash \{0\}$ (strictly positive
integer). For $t\in\Rf$ the polynomial $B[t] = B_0 + B_1 \, t +
\ldots B_d \, t^d\in\Psi^{\b}(E)$ is a classical $\pdo$ of order
$\b = \max\{\b_j \ | \ j=0,\ldots,d\}$. Suppose $A$ admits a
principal angle. Then $\zeta(A + B[t])+ \h_0(A + B[t])$ is a
polynomial of degree $d\,\left[n/(\a - \b)\right]$ in $t$ with
local coefficients. Specifically, let $Q\in\Psi^{-\a}(E)$ be a
two-sided parametrix for $A$, then

\begin{eqnarray}
\zeta(A + B[t],0)|^{{\rm mer}} + \h_0(A + B[t])  =
 \zeta(A,0)|^{{\rm mer}} + \h_0(A) \nonumber \\
   -\frac{1}{\a}\sum_{k=1}^{\left[\frac{n}{\a - \b}\right]}
\sum_{I_k} \frac{(-1)^k}{k}\ \res(Q B_{i_1} Q B_{i_2}\ldots Q
B_{i_k}) \ t^{|I_k|} \ , \label{e:zetapolynomial1}
\end{eqnarray}

\noi where the inner sum is over $k$-tuples $I_k = (i_1, \ldots,
i_k)$ of $k$ (not necessarily distinct) elements $i_j\in\{1,
\ldots, d\}$, and $|I_k| = i_1 + \ldots + i_k$. If $A$ is
invertible then $Q$ in \eqref{e:zetapolynomial1} can be replaced
by $A\ii$. In particular, if $A$ is invertible
\begin{equation}\label{e:zetapolynomial2}
\zeta(A + tI,0)|^{{\rm mer}} = \zeta(A,0)|^{{\rm mer}}
-\frac{1}{\a}\sum_{k=1}^n \frac{(-1)^k}{k} \ \res(A^{-k}) \ t^k \
.
\end{equation}
\end{thm}

\vskip 3mm Note that \eqref{e:zetapolynomial1} implies that if
$B\in\Psi^{\b}(E)$ and $\a - \b - n \in \Nf\backslash \{0\}$, then
$$\zeta(A + B,0)|^{{\rm mer}} + \h_0(A + B) = \zeta(A,0)|^{{\rm
mer}} + \h_0(A) \ .$$ Equation \eqref{e:zetapolynomial2} is
equivalent to equation (6.5) of \cite{Wo1}.

\vskip 5mm

\noi {\bf Example}. \ \ If we apply \eqref{e:zetapolynomial2} to a
Laplace-type differential operator $\Delta_g$ on a
$2m$-dimensional manifold, then comparing with \eqref{e:zetazerot}
we infer for $k\geq 1$
\begin{equation}\label{e:c0res}
\frac{1}{2} \ \res(\Delta_g^{-k}) =
\frac{c_{-k}(\Delta_g)}{(k-1)!} \ .
\end{equation}
Specifically, for $k=m$
\begin{equation*}
\frac{1}{2} \ \res(\Delta_g^{-m}) = \frac{\Vol_g(M)\, {\rm
rk}(E)}{(4\pi)^m \, (m-1)!}
\end{equation*}
which coincides with \eqref{e:resdeltam}, while, for example, for
the Laplace-Beltrami operator on a $4$-manifold with scalar
curvature $\k_M$
\begin{equation*}
\ \res(\Delta_g^{-2}) = \frac{1}{24\pi^2}\int_M \k_M \ dx \ .
\end{equation*}
\vskip 1mm

With a little extra work it is easy to see using these methods for
$A\in\Psi^{\a\neq 0}(E)$ not necessarily positive but admitting a
principal angle, that

$$ \frac{1}{\a} \ \res(A^{-z_0}) = {\rm Res}^{\Cf}_{s=z_0}
\z(A,s)|^{{\rm mer}}_{s=z_0} \ , $$

\vskip 3mm

\noi where the right side is the usual complex residue of a
meromorphic function. This formula is well-known \cite{Wo2} and
\eqref{e:c0res} is a particular case.

%*****************************************************

\vskip 5mm

\section{Preliminaries}

\vskip 3mm To explain the nature of the residue determinant we
begin with the sub-algebra $\Psi^0(E)$ of classical $\pdo$s of
order zero. For $A\in\Psi^0(E)$ with spectrum disjoint from the
ray $R_{\th} = \{re^{i\th} \ | \ r\geq 0\}$ the complex powers for
any $s\in \Cf$ are defined by
\begin{equation}\label{e:0powers}
A_{\th}\si = \frac{i}{2\pi}\int_{\Gamma_{\theta}}\la_{\th}\si
(A-\la I)\ii \ d\la \ \ \in \ \ \Psi^0(E) \ ,
\end{equation}
where $\Gamma_{\theta}$ is a bounded `keyhole' contour enclosing
the spectrum of $A$ but not enclosing any part of $R_{\th}$ (or
the origin), while the logarithm of $A$ is defined by
\begin{equation}\label{e:0log}
    \log_{\th}A = \frac{i}{2\pi}\int_{\Gamma_{\theta}}\log_{\th}\la \
    (A-\la I)\ii \ d\la \ \  \in \ \ \Psi^0(E) \ .
\end{equation}
Here $\la_{\th}\si, \ \log_{\th}\la := -(d/ds)\la_{\th}\si|_{s=0}$
are the branches defined by $\la_{\th}\si = |\la|\si
e^{-is\arg(\la)}$ with $\th - 2\pi \leq \arg(\la) < \th$. These
formulas are valid in any Banach algebra.

Sitting inside $\Psi^0(E)$ is the ideal $\Psi^{<-n}(E)$. An
operator $Q\in\Psi^{<-n}(E)$ is trace class with an absolutely
integrable matrix-valued kernel
\begin{equation}\label{e:kernel}
K(Q,x,y) = \frac{1}{(2\pi)^n}\int_{\Rf^n} e^{i(x-y).\xi}
    \ \s(Q)(x,\xi) \,d\xi
\end{equation}
over $M\times M$ smooth away from the diagonal (so for all
$\pdo$s) and continuous along $\{(x,x) \ | \ x\in M\}$.
Consequently, $Q$ has $L^2$ trace
\begin{equation}\label{e:L2trace}
    \Tr(Q) = \int_M \tr(\,K(Q,x,x)\,) \, dx \ .
\end{equation}
This is a {\em non-local} spectral invariant; from
\eqref{e:kernel} the trace \eqref{e:L2trace} depends on the
complete symbol $\s(Q)$, not just on finitely many homogeneous
terms.

Taken together, \eqref{e:0log} and \eqref{e:L2trace} define a
determinant on the ring $I + \Psi^{<-n}(E)$ of $\pdo$s which
differ from the identity by an element of $\Psi^{<-n}(E)$
\begin{equation}\label{e:TrlogQ}
 {\rm det}_{{\rm {\Small Tr}}} : I + \Psi^{<-n}(E)
\too \Cf, \hskip 5mm
 \log{\rm det}_{{\rm {\Small Tr}}}(A) := \Tr(\log_{\th}(A)) \ .
 \end{equation}
Here, one uses
\begin{equation}\label{e:logQtrclass}
A\in I + \Psi^{<-n}(E) \ \ \ \Rightarrow \ \ \ \log_{\th}A \in
\Psi^{<-n}(E) \ .
\end{equation}

${\rm det}_{{\rm {\Small Tr}}}$ extends the usual determinant in
finite-dimensions to the group $(I + \Psi^{<-n}(E))_{{\rm inv}}$
of invertible operators in $I + \Psi^{<-n}(E)$:

\vskip 2mm

\begin{lem}\label{lem:classicalres}
On the group $(I + \Psi^{<-n}(E))_{{\rm inv}}$ the determinant
${\rm det}_{{\rm {\Small Tr}}}$ is the Fredholm determinant; it is
independent of the choice of $\th$ and is multiplicative. ${\rm
det}_{{\rm {\Small Tr}}}$  is a non-local spectral invariant, and
has no (multiplicative) extension from $(I + \Psi^{<-n}(E))_{{\rm
inv}}$  to the group $\Psi^{*}(E)$ of invertible elliptic $\pdo$s.
\end{lem}
\begin{proof}
The first sentence follows, for example, from the comments around
\eqref{e:relzeta}, below, and \cite{Sc} Prop(2.21). Alternatively,
one can prove these properties directly; independence of $\theta$
by the method of proof of \propref{p:localnotheta},
multiplicativity by the Campbell-Hausdorff theorem (cf. proof of
\thmref{t:mult}). Since the $L^2$ trace $\Tr$ is non-local and has
no extension to a trace functional on $\Psi(E)$, the second
statement easily follows.
\end{proof}

\vskip 1mm

\noi Since $\log_{\th}A$ vanishes on $\Ker(A)$, the functional
${\rm det}_{{\rm {\Small Tr}}}$ is discontinuous at non-invertible
elements of the ring $I + \Psi^{<-n}(E)$, and therefore differs
from the Fredholm determinant which is continuous and vanishes on
such elements. On the other hand, from \eqref{e:restracezero} and
\eqref{e:logQtrclass}, or from \eqref{e:detresstable}:

\vskip 1mm

\begin{lem}
The residue determinant is trivial on $I + \Psi^{<-n}(E)$
\begin{equation}\label{e:detrestrivonfreddet}
A\in I + \Psi^{<-n}(E) \ \ \ \Rightarrow \ \ \ {\rm det}_{{\rm
\res}}A = 1 \ .
\end{equation}
\end{lem}

\vskip 1mm

The relation of the residue determinant \eqref{e:resdet} to the
classical determinant \eqref{e:TrlogQ} is thus structurally the
same as that of the residue trace to the classical $L^2$ operator
trace.

\vskip 2mm

\noi Because adding a smoothing operator to $A\in\Psi(E)$ does not
affect ${\rm det}_{{\rm res}}A$ the invertibility of the operator
is not detected in either \eqref{e:degreezero} or
\eqref{e:detrestrivonfreddet}.

\vskip 4mm

Next, let $A$ be a classical $\pdo$  of order $\a
> 0$ with principal angle $\th$. We assume further for simplicity
that $\th$ is an Agmon angle, meaning that $A - \la I$ is
invertible in a neighborhood of $R_{\th}$; in particular, $A$ is
elliptic and invertible. Since the $L^2$ operator norm of $(A-\la
I)\ii$ is $O(|\la|\ii)$ as $|\la|\to\o$ one can define for $\re(s)
> 0$
\begin{equation}\label{e:powers}
A_{\th}\si = \frac{i}{2\pi}\int_{\Cc_{\theta}}\la_{\th}\si (A-\la
I)\ii \ d\la \ , \hskip 10mm \re(s) > 0 \ ,
\end{equation}
where now $\Cc_{\theta}$ is a contour travelling in along
$R_{\th}$ from infinity to a small circle around the origin,
clockwise around the circle, and then back out along $R_{\th}$ to
infinity.

\vskip 3mm

In Seeley's 1967 paper \cite{Se} on the complex powers
$A\si_{\th}$ of an elliptic operator, two quite different
extensions of \eqref{e:powers} to the whole complex plane were
explained.

\vskip 3mm

In the first of these, for $s\in\Cf$ choosing $k\in\Nf$ with
$\re(s) + k> 0$ and setting
\begin{equation}\label{e:ext1}
A_{\th}\si :=  A^k A_{\th}^{-s-k} \ \ \in \ \ \Psi^{-\a s}(E) \ ,
\end{equation}
where $A_{\th}^{-s-k}$ on the right side is defined by
\eqref{e:powers}, defines, independently of $k$, a group of
elliptic classical $\pdo$s with
\begin{equation}\label{e:powersAth}
A_{\th}^0 = I, \ A_{\th}^m = A^m, \ m\in\Zf \ .
\end{equation}
An important consequence is the construction of the logarithm of
$A$, a logarithmic $\pdo$ of type $\a$ \cite{Ok1,KV}, defined by
\begin{equation}\label{e:logA}
\log_{\th}A :=  -\left.\frac{d}{ds}\right|_{s=0}A_{\th}\si \ ,
\end{equation}
and satisfying $\frac{d}{ds}A_{\th}\si = -\log_{\th}A \cdot
A_{\th}\si$.  Since $\log_{\th}A$ and $A_{\th}\si$ for $\re(s)<0$
are unbounded (note \eqref{e:powersAth}), and hence far from trace
class, this extension has been of little direct interest. It is,
though, precisely the object of interest, and, with the residue
trace at hand, all that is needed to define the residue
determinant.

\vskip 3mm

Traditionally, however, one does something quite different and
attempts to extend the determinant \eqref{e:TrlogQ} from $I +
\Psi^{< -n}(E)$ to $\Psi(E)$ by extending the $L^2$ trace from
$\Psi^{< -n}(E)$ to $\Psi(E)$ using spectral zeta-functions. The
latter leads necessarily to a quasi-trace on $\Psi(E)$ (i.e. a
functional which is non-tracial), and hence to a quasi-determinant
(i.e. a functional which is non-multiplicative), the
zeta-determinant. This is achieved through the meromorphic
extension of the complex powers $A_{\th}\si$ from $\re(s) > n/\a$,
the second of Seeley's extensions, and has been the object of
enormous interest. It is, though, quite irrelevant to the
construction of the residue determinant.

\vskip 3mm

Nevertheless, spectral zeta functions have a part to play in what
follows and we need to recall something of these constructions. We
include a $\pdo$ coefficient $B$, which for the moment we assume
to be classical of order $\b$, and we assume $A$ to have order
$\a>0$ and to be invertible. The kernel $K(B\,A_{\th}\si,x,y)$ of
$BA_{\th}\si$, which is continuous in $(x,y)$ and holomorphic in
$s$ for $\re(s) > (n+\b)/\a$, has along the diagonal a meromorphic
extension $K(B\, A_{\th}\si,x,x)|^{{\rm mer}}$ to all $s\in\Cf$
with at most simple poles located on the real axis at the points
indicated in \eqref{e:zetapoles}. Consequently the zeta function
$\z_{\th}(B,A,s) := \Tr(B\,A_{\th}\si)$ is holomorphic for
$\re(s)> (n+\b)/\a$ and extends to a meromorphic function
\begin{equation*}
\z_{\th}(B,A,s)|^{{\rm mer}} := \int_M
\tr(K(B\,A_{\th}\si,x,x)|^{{\rm mer}}) \ dx \ \
\end{equation*}
on $\Cf$. It has pole structure (\cite{Se}, \cite{GS}, \cite{Gr})
\begin{equation}\label{e:zetapoles}
\G(s)\,\z_{\th}(B,A,s)|^{{\rm mer}} \sim
\sum_{j=0}^{\o}\frac{c_j}{ s + \frac{j - n - \b}{\a}}  +
\sum_{k=0}^{\o}\left(\frac{c^{\prime}_k}{(s+ k)^2} +
\frac{c^{\prime\prime}_k}{s+ k}\right) \ ,
\end{equation}
where the terms $c_j, c^{\prime}_j$, are local, depending on just
finitely many homogeneous terms of the symbols of both $A$ and
$B$, while the $c^{\prime\prime}_k$ are global, depending on the
complete symbols. Around zero \eqref{e:zetapoles} implies a
Laurent expansion $$ \,\z_{\th}(B,A,s)|^{{\rm mer}}  =
\frac{c^{\prime}_0}{s}  + ( c_{n + \b} +  c^{\prime\prime}_0) +
O(s) \ , \hskip 10mm {\text as} \ \ s\too 0 \ ,$$ the simple pole
determining the zeta-function formula for the residue trace of the
classical $\pdo$ $B$
\begin{equation}\label{e:zetares}
\res (B)  = \a \, {\rm Res}^{\Cf}_{s=0} \, (\,
\z_{\th}(B,A,s)|^{{\rm mer}}\, )= \a\, c^{\prime}_0 \ ,
\end{equation}
while  the constant term defines the `$A$-weighted zeta-trace'
\begin{equation}\label{e:zetatrace}
\Tr_{\zeta}^A(B) = c_{n + \b} +  c^{\prime\prime}_0 \ .
\end{equation}
If $\b< - n$ one has $\res(B) = 0$ and so $\z_{\th}(B,A,s)|^{{\rm
mer}}$ is then holomorphic near $s=0$, while  $\Tr_{\zeta}^A(B) =
\z_{\th}(B,A,0)|^{{\rm mer}}$ and is equal to the $L^2$ trace
$\Tr(B) = c^{\prime\prime}_0$. If $\b\in \Rf\backslash \Zf$ then
once more $\res(B) = 0$ and $\Tr_{\zeta}^A(B)
=\z_{\th}(B,A,0)|^{{\rm mer}}$ is independent of $A$ \cite{KV},
and again equal to the global term $c^{\prime\prime}_0$ \cite{Gr},
and vanishes on commutators for which the sum of operator orders
is non-integral \cite{KV}. These properties also hold on the
subalgebra of odd-class $\pdo$s \cite{KV}. $\Tr_{\zeta}^A$ is a
quasi-trace in so far as it is not tracial on the full algebra
$\Psi(E)$, but does extend the $L^2$ trace to a trace on the above
subclasses.

\vskip 2mm

If $B$ is a {\em logarithmic} $\pdo$ then $\z_{\th}(B,A,s)|^{{\rm
mer}}$ again extends meromorphically to $\Cf$ but now with $\b =
0$ in \eqref{e:zetapoles} and with possible additional poles
$c_{j,1}(s+ \frac{j-n}{\a})^{-2}$ \cite{Gr}. When $B =
\log_{\th}A$ there is no pole at $s=0$ and a (quasi-) determinant,
the zeta determinant ${\rm det}_{\z,\th}(A)$, may be defined by
taking the zeta trace \footnote{For $L$ a logarithmic $\pdo$, such
as $\log_{\th}A$, the dependence on the choice of regularizing
operator $\Tr_{\zeta}^{A_1}(L) - \Tr_{\zeta}^{A_2}(L) =
\Tr(L(A_1\si-A_2\si))|_{s=0}^{{\rm mer}}$  is computed in
\cite{KV}, \cite{Ok2}, \cite{PaSc} as a residue trace, leading to
the multiplicative anomaly formula for the zeta-determinant.} of
$\log_{\th}A$
\begin{eqnarray}
\log{\rm det}_{\z,\th}(A) & := & \Tr_{\zeta}^A(\log A)  \nonumber \\
& = & \z_{\th}(\log A, A,0)|^{{\rm mer}} \ . \label{e:zetadet1}
\end{eqnarray}

\vskip 2mm

\noi If $\a >0$ and one sets $\z(A,s)|^{{\rm mer}} :=
\z(I,A,s)|^{{\rm mer}} =  \Tr(A\si)|^{{\rm mer}}$, then
equivalently
\begin{equation}\label{e:zetadet2}
\log{\rm det}_{\z,\th}(A)  =  -
\frac{d}{ds}\z_{\th}(A,s)|_{s=0}^{{\rm mer}} \ , \hskip 10mm \a >
0 \ .
\end{equation}
\noi If $\a = 0$ then  from \eqref{e:logQtrclass} and the comment
following \eqref{e:zetatrace}, the determinant \eqref{e:zetadet1}
is defined on $I + \Psi^{<-n}(E)$ and equal there to the Fredholm
determinant.  On the other hand, $\z_{\th}(A,s)$, and hence the
right side of \eqref{e:zetadet2}, is then not defined for any $s$.
Nevertheless, the relative zeta function
\begin{equation}\label{e:relzeta}
\z_{\th}^{{\rm rel}}[A,B](s) := \Tr(A_{\th}\si - B_{\th}\si)
\end{equation}
is defined on $I + \Psi^{<-n}(E)$, and one has there
$\Tr\log_{\th}(A) = -\frac{d}{ds}\z_{\th}^{{\rm
rel}}[A,I](s)|_{s=0}$. For elliptic $A_0, A_1$ of non-zero order
$\z_{\th}^{{\rm rel}}[A_0,A_1](s) = \z_{\th}(A_0,s) -
\z_{\th}(A_1,s)$; relative determinants are studied in \cite{Mu},
\cite{Sc}.

\vskip 3mm

The zeta determinant ${\rm det}_{\z,\th}(A)$  is non-local,
depends on the spectral cut $R_{\th}$ \cite{Wo1} and is not
multiplicative \cite{Ok2,KV}.

\vskip 3mm

The residue determinant has a quite different nature.

\vskip 3mm

First, the {\em construction} of the residue determinant takes
place completely independently of spectral zeta-functions, and it
is in this distinction that its non-triviality lies.

\vskip 2mm

Specifically, using the spectral $\z$-function formula
\eqref{e:zetares} with $B = \log_{\th}A$ to define a putative
residue determinant, rather than the symbolic definition
\eqref{e:resdet}, leads to a trivial determinant (equal to $1$);
the triviality is equivalent to $\z(A,s)|^{{\rm mer}}$ being
holomorphic at zero, and thus \eqref{e:zetadet1} being defined.

Further, since the residue zeta function $\z_{\res}(A,s) :=
\res(A\si)$ is highly discontinuous -- for, from
\eqref{e:restracezero}, $\res(A\si)$ can be non-zero only for $s
\cdot\ord(A) \in\Zf\cap (\,-\o \ , \, n\,]$ -- there is no residue
analogue of \eqref{e:zetadet2}.

\vskip 3mm

The residue trace \eqref{e:resdet} on $\log_{\th}A$ thus does not
arise a complex residue, but it does generalize the integral
formula which for classical $\pdo$s coincides with the complex
residue \eqref{e:zetares}. However, if $A, B$ have residue
determinants and are of the same order then the difference
$\res(\log A) -\res(\log B)$ is given as a complex residue.

\vskip 2mm

The residue  determinant ${\rm det}_{\res}(A)$  is local,
independent of the spectral cut $R_{\th}$ and multiplicative.

%*****************************************************

\vskip 5mm

\section{Proofs}

\vskip 3mm

Let $\as \sim (\as_0, \as_1, \ldots)\in S^{\a}(U), \ \bs \sim
(\bs_0, \bs_1, \ldots)\in S^{\b}(U)$ be local classical (1-step
polyhomogeneous) symbols of respective degrees $\a, \b\in\Rf$.
Then a product structure is defined on $S(U)$
$$\as \circ \bs \sim ((\as \circ \bs)_0, (\as \circ \bs)_1, \ldots \ )\in
S^{\a+\b}(U) \ ,$$ with
\begin{equation}\label{e:symbolproduct}
(\as \circ \bs)_j = \sum_{|\mu| + k + l = j} \frac{1}{\mu
!}\partial_{\xi}^{\mu}(\as_k) \, D_x^{\mu}(\bs_l) \ ,
\end{equation}
and multiplicative identity element
$$\Is = (I,0,0,\ldots) \ .$$

At the $\pdo$ level this represents the operator product modulo
smoothing operators; thus if in local coordinates $\s(A), \s(B)$
are respectively equivalent symbols to $\as, \bs$, then $\s(AB)$
is equivalent to $\as\circ\bs$. This is all that is needed to
compute local quantities such as the residue determinant.

To do so for a classical elliptic $\pdo$ $A$ of order $\a$,
standard methods \cite{Gi,Se,Sh} construct a parametrix for $A
-\la I$ by inverting locally at the symbolic level. We consider a
finite open cover of $M$ by coordinate patches $U_i, i\in I =
\{1,\ldots,m\}$, over each of which $E$ is trivialized, with a
subordinate partition of unity $\phi_i\in\Ci(U_i)$ such that for
$i,j\in I$ there is a $l_{ij}\in I$ with $\supp(\phi_i) \cup
\supp(\phi_j) \subset U := U_{l_{ij}}$. Then
\begin{equation}\label{e:localpdo}
A = \sum_{i,j}\phi_i A \phi_j
\end{equation}
with each summand a $\pdo$ acting in a single coordinate patch,
and it will be enough to work with a symbol $ \as = \s(\phi_i A
\phi_j)\sim (\as_0, \as_1, \ldots)  \in S^{\a}(U)$ of each such
local operator. A local resolvent symbol
$$\rs(\la) \sim (\rs(\la)_0, \rs(\la)_1, \ldots \ )\in S^{-\a}(U_{\la}) $$
is defined over $U_{\la} = \{x\in U \ | \ \la \notin {\rm
spec}(\as_0(x,\xi)), \ ,\xi\in\Rf^n\}$ by the inductive formulae
\begin{equation}
\rs(\la)_0(x,\xi) = (\as_0(x,\xi) - \la \Is)\ii \ ,
\end{equation}
\begin{equation}\label{e:symbolinverse}
\rs(\la)_j(x,\xi) =  -\,\rs(\la)_0(x,\xi)\sum_{\stackrel{|\mu| + k
+ l = j}{l<j}} \frac{1}{\mu !}\partial_{\xi}^{\mu}\as_k(x,\xi) \,
D_x^{\mu}\rs(\la)_l(x,\xi)  \ .
\end{equation}
For $|\xi|\geq 1, t\geq 1$, each $\rs(\la)_j(x,\xi)$ has the
quasi-homogeneity property
\begin{equation*}
\rs(t^{\a}\la)_j\,(x,t\,\xi) = t^{-\a-j}\,\rs(\la)_j(x,\xi) \ ,
\end{equation*}
and by construction
\begin{equation}\label{e:symbolinverse2}
\rs(\la)\circ (\as - \la \Is) \sim  (\as - \la \Is)\circ \rs(\la)
\sim \Is \ .
\end{equation}
Consequently, if $A$ has principal angle $\th$, then $\log_{\th}A$
and $A_{\th}\si$ are represented by symbols
\begin{equation}\label{e:logpowersymbol}
    \log_{\th}\as \sim (\qs_{\th,0}, \qs_{\th,1}, \ldots \ ) \ , \hskip
    10mm
    \as_{\th}\si \sim (\as_{\th,0}\si, \as_{\th,1}\si, \ldots \ ) \ ,
\end{equation}
where, with $\Gamma_{\theta} = \Gamma_{\theta}(x,\xi)$ a closed
contour as in \eqref{e:0powers} chosen to enclose the spectrum of
$\as_0(x,\xi)$ avoiding the spectral cut and the origin,
\begin{equation}\label{e:logsymbol}
\qs_{\th,j}(x,\xi) =
\frac{i}{2\pi}\int_{\Gamma_{\theta}}\log_{\th}\la \
\rs(\la)_j(x,\xi) \ d\la \ ,
\end{equation}
\begin{equation}\label{e:powersymbol}
\as_{\th,j}\si(x,\xi) =
\frac{i}{2\pi}\int_{\Gamma_{\theta}}\la_{\th}\si \
\rs(\la)_j(x,\xi) \ d\la \ .
\end{equation}
We may also write
$$ (\log \as)_{\th,j}(x,\xi) := \qs_{\th,j}(x,\xi) \ .$$
For $|\xi|\geq 1$, it follows for $j\geq 0$ that
$\as_{\th,j}\si\in S^{-\a s - j}(U)$ is homogeneous of degree $-\a
s - j$, and for $j\geq 1$ that $\qs_{\th,j}\in S^{- j}(U)$ is
homogeneous of degree $- j$. Since
\begin{equation*}
    \qs_{\th,j}(x,\xi)  = -\partial_{s}|_{s=0}\as_{\th,j}\si(x,\xi)
\end{equation*}
and
\begin{equation*}
\as_{\th,j}\si(x,\xi)|_{s=0} = \d_{j,0}\Is \ ,
\end{equation*}
then
\begin{equation}\label{e:symbolpowerzero}
\as_{\th}^0  = \Is = (I,0,0,\ldots)
\end{equation}
and
\begin{equation}\label{e:logsymbolzero}
\qs_{\th,0}(x,\xi) = \a\log [\xi] I + \ps_{\th,0}(x,\xi)
\end{equation}
with $\ps_{\th,0}\in S^0(U)$ a classical symbol of degree $0$.
This means that the $\pdo$s $\log_{\th}A$ and $A_{\th}\si$ can be
approximated modulo smoothing operators as
\begin{equation}\label{E:approxmodsmoothing}
\log_{\th}A \sim \sum_{j=0}^{\o} (\log_{\th}A)_{[j]} \ , \hskip
10mm A_{\th}\si \sim \sum_{j=0}^{\o} A_{\th,j}^{(-s)}
\end{equation}
with $(\log_{\th}A)_{[j]} = {\rm OP}(\qs_{\th,j})$, an operator in
$\Psi^{-j}(U, E)$ for $j\geq 1$, and $A_{\th,j}^{(-s)} = {\rm
OP}(\as_{\th,j}\si)\in \Psi^{-s\a - j}(U, E)$ for $j\geq 0$. In
particular, the local residue density associated to
$\qs_{\th,n}(x,\xi)$ is defined independently of the choice of
local coordinates and
\begin{equation}\label{e:resdet3}
\log{\rm det}_{\res} A  :=
\frac{1}{(2\pi)^n}\int_M\int_{|\xi|=1}\tr(\qs_{\th,n}(x,\xi)) \ d
S(\xi)\, dx \ .
\end{equation}

\vskip 5mm

\noi {\bf Proof of \propref{p:localnotheta}}

\vskip 1mm

\begin{proof}
From \eqref{e:symbolinverse} we have that $\rs(\la)_n$ is computed
only using the first $n+1$ homogeneous terms $\as_0, \ldots,
\as_n$. Consequently, by \eqref{e:logsymbol} and \eqref{e:resdet3}
the same is true for the residue determinant.

\vskip 2mm

Let $\th, \phi\in\Rf$ be two choices of principal angle with $(\th
- \phi)/2\pi\in\Rf\backslash\Zf$. Then
\begin{equation}\label{e:qthphi}
\qs_{\phi,j}(x,\xi)  - \qs_{\th,j}(x,\xi) = 2\pi i\, m \,\Is_j  \
+ \ \int_{\G_{\phi,\th}} \rs(\la)_j(x,\xi) \ d\la \ ,
\end{equation}
where $m= \pm[(\th - \phi)/2\pi]\in\Zf$ and the bounded contour
$\G_{\phi,\th} = \G_{\phi,\th}(x,\xi)$ can be taken of the form
$$ \{\rho e^{i\th} \ | \ R \geq \rho \geq r \} \cup \{ r\,e^{it} \
| \ \phi \geq t \geq \th \} \cup \{\rho e^{i(\th- 2\pi)} \ | \ r
\leq \rho \leq R \} \cup \{ R\,e^{it} \ | \ \phi \leq t \leq \th
\} $$  enclosing an annular region between the cuts $R_{\th}$, and
$R_{\phi}$ and circles of radius $r<R$. This follows by a similar
analysis to \cite{Wo1}\S3 for the symbols of the complex powers.

On the other hand, the contour integral on the right-side of
\eqref{e:qthphi} is $-2\pi i$ times the homogeneous component of
degree $-j$ of the local symbol of a $\pdo$ projection
$P_{\th,\phi}(A)$ whose range contains the direct sum of those
generalized eigenspaces of $A$ with eigenvalues contained in
$\G_{\phi,\th}$, and is zero if $(\th - \phi)/2\pi\in\Zf$ (see
\cite{Bu}, \cite{Po}). Consequently, taking $j=n$,
\eqref{e:resdet3} and \eqref{e:qthphi} imply
\begin{equation*}\label{e:thphi}
{\rm res} (\log_{\th}A) -  {\rm res} (\log_{\phi}A)  = - 2\pi
i\,{\rm res} (P_{\th,\phi}(A)) \ .
\end{equation*}
Since the residue trace of any $\pdo$ projection is zero
\cite{Wo2}, we infer that $\det_{\rm res}$ is independent of the
choice of principal angle.
\end{proof}

\vskip 1mm

\begin{rem}
The vanishing of $\res$ on $\pdo$ projections is shown in
\cite{Wo2} to be equivalent to $\z_{\th}(A,0)|^{{\rm mer}}$ being
independent of $\th$ .
\end{rem}

\vskip 10mm

\noi {\bf Proof of \thmref{t:mult}}

\vskip 2mm

\begin{proof} From \eqref{e:resdet3} and
\eqref{e:symbolproduct}, $\res(\log A)$ is seen to depend on only
the first $n+1$ homogeneous terms in the local symbol expansion of
$A$, and finitely many of their derivatives, while
$(\as\circ\bs)_n$ is determined using only $\as_0,\ldots ,\as_n,
\bs_0,\ldots ,\bs_n$. The demonstration of multiplicativity can
therefore be reduced to a certain finite-dimensional symbol
algebra, introduced by Okikiolu \cite{Ok1}\S3, where the following
standard Banach algebra version of the Campbell-Hausdorff Theorem
\cite{Ok1,J} can be applied.

\vskip 2mm

 \noi {\bf Theorem} \  {\it Let $\Bb$ be a Banach
algebra with norm $\| \ . \ \|$ and identity $I$. For invertible
elements $a, b\in\Bb$ and a choice of Agmon angles one can define
using \eqref{e:0log} elements $\log (a), \,\log (b)$ and $\log
(ab)$ in $\Bb$. Then for real sufficiently small $s,t>0$
\begin{equation}\label{e:banachcampbellhausdorff}
\log(a^s \, b^t) = s\,\log(a) + t\,\log(b) +
\sum_{k=1}^{\o}C^{(k)}(s\,\log(a),t\,\log(b)) \ ,
\end{equation}
where $C^{(k)}(s\,\log(a),t\,\log(b))$ is the element of $\Bb$
\begin{equation}\label{e:Ck}
\sum_{j=1}^{\o}\frac{(-1)^{j+1}}{j+1}\sum \frac{({\rm
Ad}(s\,\log(a))^{n_1} ({\rm Ad}(t\,\log(b))^{m_1}\ldots ({\rm
Ad}(s\,\log(a))^{n_j} ({\rm
Ad}(t\,\log(b))^{m_j}\log(b)}{(1+\sum_{i=1}^j m_i)\, n_1 ! \ldots
n_j !\,  m_1 ! \ldots m_j !}
\end{equation}
and the inner sum is over $j$-tuples of pairs $(n_i,m_i)$ such
that $n_i + m_i > 0$ and $\sum_{i=1}^j n_i + m_i = k$. For $c\in
\Bb$ the operator ${\rm Ad}(c)$ acts by ${\rm Ad}(c)(c^{'}) =
[c,c^{'}]$.}

\vskip 4mm

We denote by $\Ss_{[n]}(U)$ the algebra of finite-symbol sequences
of length $n$, introduced in \cite{Ok1}. An element of
$\Ss_{[n]}(U)$ is an $(n+1)$-tuple $\ps = (\ps_0,\ldots ,\ps_n)$
of polynomials
\begin{equation}\label{e:polynomials}
\ps_j : U \times \Rf^n \too \End(\Rf^N) \ , \hskip 10mm
\ps_j(x,\xi) = \sum_{|\mu| + |\nu| \leq n-1} p_{j,\mu,\nu}\,
x^{\mu}\xi^{\nu} \ ,
\end{equation}
with  $p_{j,\mu,\nu}\in \End(\Rf^N)$. $\Ss_{[n]}(U)$ is a
finite-dimensional vector space which relative to a fixed point
$(x_0,\xi_0)\in U\times\Rf^n$ can be endowed with an associative
product structure, defined for $\ps, \widetilde{\ps}\in
\Ss_{[n]}(U)$ by
\begin{equation}\label{e:nprod}
(\ps \circ \widetilde{\ps})_j  = \pi_{n-j}\left(\sum_{|\mu| + k +
l = j} \frac{1}{\mu !}\partial_{\xi}^{\mu}(\ps_k) \,
D_x^{\mu}(\widetilde{\ps}_l)\right) \ ,
\end{equation}
where for a smooth function $f$ defined in a neighborhood of
$(x_0,\xi_0)\in U\times\Rf^n$
\begin{equation}\label{e:finiteTaylor}
\pi_m(f) = \sum_{|\mu| + |\nu| \leq m}\frac{1}{\mu
!\nu!}\partial_{\xi}^{\mu}\partial^{\nu}(f)(x_0,\xi_0)\,(x-x_0)^{\mu}\,(\xi
-\xi_0)^{\nu}
\end{equation}
is the Taylor expansion of $f$ around $(x_0,\xi_0)$ to order $m$.
Endowed with this product, relative to $(x_0,\xi_0)$,
$\Ss_{[n]}(U)$ becomes an algebra which we denote by
$$\Ss_{[n]}(U)(x_0,\xi_0) \ .$$

\vskip 1mm

The map from the symbol space $S(U)$ to symbols of length $n$
\begin{equation}\label{e:StoSn}
\pi: S(U) \too \Ss_{[n]}(U)(x_0,\xi_0)\ , \hskip 5mm \pi(\as) :=
(\pi_n(\as_0), \pi_{n-1}(\as_1), \ldots , \pi_0(\as_n))
 \ ,
\end{equation}
where $\as = (\as_0, \as_1, \ldots)$, is an algebra homomorphism,
so that
\begin{equation}\label{e:StoSnhom}
 \left(\pi(\as \circ \bs)\right)_j = \left(\pi(\as) \circ
 \pi(\bs)\right)_j\ ,
\end{equation}
while, from \eqref{e:finiteTaylor}, evaluation at the point
$(x_0,\xi_0)$ gives
\begin{equation}\label{e:goodforjlessthann}
\left(\pi(\as)\right)_j(x_0,\xi_0) = \as_j(x_0,\xi_0) \ ,\hskip
5mm  j \leq n \ .
\end{equation}

\vskip 1mm

The logarithm of an element $\ps = (\ps_0,\ldots ,\ps_n)\in
\Ss_{[n]}(U)(x_0,\xi_0)$ admitting a principal angle can be
defined by the procedure used in $S(U)$: Consider, by inclusion,
$\ps$ as an element $\tilde{\ps}$ of $S(U)$. If $\la \notin {\rm
spec}(\ps_0(x,\xi))$ then $\tilde{\ps}$ has a resolvent
$\rs(\la)\in S(U_{\la}) $ given by \eqref{e:symbolinverse}, while
\begin{equation}\label{e:finitesymbolresolvent}
\rs_{\pi}(\la) := \pi(\rs(\la)) \in \Ss_{[n]}(U)(x_0,\xi_0)
\end{equation}
inverts $\ps - \la \Is_n$ in $\Ss_{[n]}(U)(x_0,\xi_0)$; that is,
since $\pi(\tilde{\ps}) = \ps$, applying $\pi$ to
\eqref{e:symbolinverse2} and using \eqref{e:StoSnhom} we have
\begin{equation}\label{e:finitesymbolinverse}
\rs_{\pi}(\la)\circ (\ps - \la \Is_n) =  (\ps - \la \Is_n)\circ
\rs_{\pi}(\la) = \Is_n \ ,
\end{equation}
where $\Is_n = (I,0,\ldots,0)$ is the identity symbol in
$\Ss_{[n]}(U)(x_0,\xi_0)$. Set
\begin{equation}\label{e:finitelogsymbol}
(\log_{\th}\ps)_{j}(x,\xi) =
\frac{i}{2\pi}\int_{\Gamma_{\theta}}\log_{\th}\la \
\rs_{\pi}(\la)_j(x,\xi) \ d\la \ .
\end{equation}
Since the entries  of $\rs_{\pi}(\la)$ are finite Taylor
expansions of $\rs(\la)_j(x,\xi)$ around $(x_0,\xi_0)$ the only
logarithmic term is a $\log|\xi_0|$, there is no $\log|\xi|$ term.
It follows that \eqref{e:finitelogsymbol} is an element of
$\Ss_{[n]}(U)(x_0,\xi_0)$. Moreover, it is clear (\cite{Ok1} Lemma
3.6) that for $\as\in S(U)$
\begin{equation}\label{e:pilog=logpi}
\left(\pi(\log_{\th}\as)\right)_j =
\left(\log_{\th}(\pi(\as)\right)_j \ .
\end{equation}
Likewise, $\ps_{\th,j}\si(x,\xi) =
\frac{i}{2\pi}\int_{\Gamma_{\theta}}\la\si_{\th}
\rs_{\pi}(\la)_j(x,\xi) \ d\la \in\Ss_{[n]}(U)(x_0,\xi_0) $ if
$\ps\in\Ss_{[n]}(U)(x_0,\xi_0)$, and we find for $\as\in S(U)$
\begin{equation}\label{e:pipower=powerpi}
\left(\pi(\as_{\th}\si)\right)_j =
\left(\,(\pi(\as))_{\th}\si\,\right)_j \ .
\end{equation}

\vskip 1mm

Now for $s,t\in [0,1]$ and $\as, \bs\in S(U)$, \eqref{e:StoSnhom},
\eqref{e:pilog=logpi} and \eqref{e:pipower=powerpi} give
\begin{equation}\label{e:pilogprod=logpiprod}
\left(\,\pi(\log_{\th}(\as^s\circ\bs^t))\,\right)_n   =
\left(\,\log_{\th}(\pi(\as)^s\circ\pi(\bs)^t)\,\right)_n \ ,
\end{equation}
omitting the $\th$ subscript. Since $\Ss_{[n]}(U)(x_0,\xi_0)$ is a
finite-dimensional algebra, for $s,t \geq 0$ sufficiently small we
have from \eqref{e:banachcampbellhausdorff} for the induced norm
\begin{eqnarray}
\left(\,\pi(\log_{\th}(\as^s\circ\bs^t))\,\right)_n  & = &
s\,(\log\,\pi(\as))_n + t\,(\log\,\pi(\bs))_n \\ & &  +
\sum_{k=1}^{\o}\left(C^{(k)}(s\,\log\,\pi(\as), \
t\,\log\pi(\bs))\right)_n  \nonumber \\
& = & s\,(\pi(\log\,\as))_n + t\,(\pi(\log\,\bs))_n \nonumber\\ &
& + \sum_{k=1}^{\o}\left(\pi\left(C^{(k)}(s\,\log\,\as, \ t\,\log
\,\bs)\right)\,\right)_n \ , \nonumber
\end{eqnarray}
and so, evaluating at the point $(x_0,\xi_0)$,
\eqref{e:goodforjlessthann} implies
\begin{eqnarray}
\log_{\th}(\as^s\circ\bs^t)_n(x_0,\xi_0) & = &
s\,\log\,\as_n(x_0,\xi_0) + t\,\log\,\bs_n(x_0,\xi_0) \label{e:x0xi0}\\
& + & \sum_{k=1}^{\o}\left(C^{(k)}(s\,\log\,\as, \ t\,\log
\,\bs)\right)_n(x_0,\xi_0) \ . \nonumber
\end{eqnarray}
Since all terms in \eqref{e:x0xi0} lie in the symbol class $S(U)$,
with uniformly continuous derivatives of all orders on compact
subsets of $U\times \Rf^n$, the convergence in \eqref{e:x0xi0} as
$N\to \o$ of

$$\log_{\th}(\as^s\circ\bs^t)_n -
s\,\log\,\as_n - t\,\log\,\bs_n -
\sum_{k=1}^{N}\left(C^{(k)}(s\,\log\,\as, \ t\,\log
\,\bs)\right)_n$$

\noi and all its derivatives at the point $(x_0,\xi_0)$ is also
uniform in $(x_0,\xi_0)\in U_c\times S^{n-1}$ for compact subsets
$U_c \subset U$. Hence, taking a partition of unity we can
interchange the sum with integration over $S^*M$ to get

\begin{equation}\label{e:ressum}
\int_M\int_{|\xi|=1}\tr(\log(\as^s\circ\bs^t)_n(x,\xi)) \ d
S(\xi)\,dx =
\end{equation}
$$ s\int_M\int_{|\xi|=1}\tr(\log\,\as_n(x,\xi)) \ d S(\xi)\,dx
+ t\int_M\int_{|\xi|=1}\tr(\log\,\bs_n(x,\xi)) \ d S(\xi)\,dx $$
$$+ \sum_{k=1}^{\o}\int_M\int_{|\xi|=1}\tr(C^{(k)}(s\,\log\,\as, \
t\,\log \,\bs)_n(x,\xi)) \ d S(\xi)\,dx\ .$$

\noi But $C^{(k)}(s\,\log A, \ t\,\log B )$ is  classical $\pdo$
of order $0$ with symbol
$$\s(C^{(k)}(s\,\log A, \
t\,\log B )) \sim C^{(k)}(s\,\log\,\as, \ t\,\log \,\bs) \ ,$$ and
so, in particular,
$$\s(C^{(k)}(s\,\log A, \
t\,\log B ))_n = C^{(k)}(s\,\log\,\as, \ t\,\log \,\bs)_n \ . $$
It is, furthermore, by definition a commutator of logarithmic
 $\pdo$s, and hence by \propref{p:restrace}
$$\frac{1}{(2\pi)^n}\int_M\int_{|\xi|=1}\tr(C^{(k)}(s\,\log\,\as, \
t\,\log \,\bs)_n(x,\xi)) \ d S(\xi)\,dx =
\res\left(C^{(k)}(s\,\log A, \ t\,\log B)\,\right) =0 \ .$$  Thus
\eqref{e:ressum} says that for sufficiently small $s,t\in[0,1]$
\begin{equation}\label{e:resmultst}
\res(\log(A^s B^t))   = s \, \res(\log A) + t \, \res(\log B) \ .
\end{equation}
But \eqref{e:resmultst} is analytic in $s, t$ and so it holds for
all $s,t\in [0,1]$. Evaluating at $s=t=1$ completes the proof.
\end{proof}

\vskip 5mm

\noi {\bf Proof of \thmref{t:resdetzetazero}}

\vskip 5mm

\begin{proof}
Let $\as(x,\xi) = \s(A)(x,\xi)\in S^{\a}(U)$ be the symbol of $A$
localized over $U$, as above. The complex powers
$A\si_{\th}\in\Psi^{-\a s}(E)$ are classical $\pdo$s defined in
the half-plane $\re(s)
> 0$ by \eqref{e:powers}, and elsewhere by \eqref{e:ext1}, with
local symbol
\begin{equation}\label{e:powersymbolexpansion}
\s(A\si_{\th})(x,\xi) \ \sim \ \sum_{j\geq0}\as_{\th,j}\si(x,\xi)
\ .
\end{equation}
If $A$ is not invertible, then for $s\neq 0$ \eqref{e:powersAth}
remains unchanged, while
\begin{equation}\label{e:powers2}
A_{\th}^0 = I - \Pi_0(A) \ ,
\end{equation}
with $\Pi_0(A)$ a (in general, non self-adjoint) projection onto
the generalized eigenspace $E_0(A)$ in the statement of
\thmref{t:resdetzetazero}.

\vskip 2mm

The symbol $\s(A\si_{\th})(x,\xi)$ is integrable in $\xi$ for
$\re(s)
> n/\a$ and, for such $s$, $K(A\si_{\th},x,x)\, dx$ defines a
$\Ci$ globally defined $n$-form on $M$ with values in $\End(E)$.

\vskip 1mm

\noi For $\re(s) > n/\a$ and any $J\in\Nf$ we have with
$\hat{d}\xi := (2\pi)^{-n} d\xi$
\begin{eqnarray}
K(A\si_{\th},x,x) & = & \int_{\Rf^n}\s(A\si_{\th})(x,\xi) \
\hat{d}\xi \nonumber \\
& = & \int_{\Rf^n}\left(\s(A\si_{\th})(x,\xi) - \sum_{j=0}^{J-1}
\as_{\th,j}\si(x,\xi)\right) \hat{d}\xi \ + \ \sum_{j=0}^{J-1}
\int_{\Rf^n} \as_{\th,j}\si(x,\xi) \ \hat{d}\xi  \ .
\label{e:slarge}
\end{eqnarray}
With $A\si_{\th}$ defined for all $s\in\Cf$ by \eqref{e:ext1}, the
difference
$$\s(A\si_{\th})(x,\xi) - \sum_{j=0}^{J-1}
\as_{\th,j}\si(x,\xi) \ \in \ S^{-\a\re(s) - J}(U)$$ is integrable
in $\xi$ for
\begin{equation}\label{e:whereJintegrable}
\re(s) > \frac{n-J}{\a} \ ,
\end{equation}
and so the first integral on the right-side of \eqref{e:slarge}
extends holomorphically to the half-plane
\eqref{e:whereJintegrable}. Hence choosing $J = n+1$ (or any $J >
n$) we can set $s=0$ in that integral to get, using
\eqref{e:symbolpowerzero} and \eqref{e:powers2},
\begin{eqnarray}
\left.\int_{\Rf^n}\s(A\si_{\th})(x,\xi) - \sum_{j=0}^{n}
\as_{\th,j}\si(x,\xi) \ \hat{d}\xi \right|^{{\rm mer}}_{s=0} & = &
\int_{\Rf^n}\left(\s(A^0_{\th})(x,\xi) - \sum_{j=0}^{n}
\as_{\th,j}^0(x,\xi)\right) \ \hat{d}\xi \nonumber \\
 & =
& \int_{\Rf^n}\left(\s(I - \Pi_0(A))(x,\xi) - \sum_{j=0}^{n}
\Is_{\th,j}(x,\xi)\right) \ \hat{d}\xi \nonumber \\
& = & \int_{\Rf^n} - \, \s(\Pi_0(A))(x,\xi) \ \hat{d}\xi \ .
\label{e:globalbit}
\end{eqnarray}

\vskip 1mm

The remaining objects of interest, then, are the local-kernels
along the diagonal
\begin{equation*}
\Ks_{j}^{-s}(x) = \int_{\Rf^n} \as_j\si(x,\xi) \ \hat{d}\xi\ .
\end{equation*}
Splitting the integral into two parts we have
\begin{equation}\label{e:twoparts}
\left.\Ks_{j}^{-s}(x)\right|^{{\rm mer}} = \left.\int_{|\xi|\leq
1} \as_j\si(x,\xi) \ \hat{d}\xi\right|^{{\rm mer}}  +
\left.\int_{|\xi|\geq 1} \as_j\si(x,\xi) \ \hat{d}\xi\right|^{{\rm
mer}} \ .
\end{equation}
We deal first with the second term on the right side of
\eqref{e:twoparts}, for which the symbol is homogeneous in
$|\xi|$, hence leading to only local poles (any $s$). Changing to
polar coordinates and using the homogeneity of $\as_j\si(x,\xi)$,
we have for $\re(s)>(n-j)/\a$
\begin{eqnarray}
\int_{|\xi|\geq 1} \as_j\si(x,\xi) \ \hat{d}\xi & = & \int_1^{\o}
r ^{-\a s-j+n-1} dr \int_{|\xi|=1} \as_j\si(x,\xi) \ \hat{d}S(\xi)  \label{e:part2extension1}\\
& = & \frac{1}{(\a s+j-n)}\int_{|\xi|=1} \as_j\si(x,\xi) \
\hat{d}S(\xi) \ . \label{e:part2extension2}
\end{eqnarray}
The meromorphic extension of the left side of
\eqref{e:part2extension1} is defined by \eqref{e:part2extension2}.
When $j\neq n$ then \eqref{e:part2extension2} is holomorphic
around $s=0$ and so from \eqref{e:symbolpowerzero}
\begin{eqnarray}
\left.\int_{|\xi|\geq 1} \as_j\si(x,\xi) \
\hat{d}\xi\right|_{s=0}^{{\rm
mer}} & = & 0 \ , \hskip 10mm j\neq 0,\, n \ , \label{e:jnotzeronorn} \\
\left.\int_{|\xi|\geq 1} \as_0\si(x,\xi) \
\hat{d}\xi\right|_{s=0}^{{\rm mer}} & =& - \,
\frac{1}{(2\pi)^n}.\frac{{\rm vol}(S^{n-1})}{n}   \ .
\label{e:jzero}
\end{eqnarray}

\vskip 2mm

For $j=n$ we use \cite{Ok2}\,Lemma(2.1) which states that there is
an equality
\begin{equation*}
\as_j\si(x,\xi) = \sum_{k=0}^{\o}
\frac{(-s)^k}{k!}((\,\log\as)^k)_j(x,\xi) \ ,
\end{equation*}
where $(\log\as)^k := \log\as\circ \log\as\circ\ldots\circ\log\as$
($k$ times) and  the right-side is convergent as a function of
$(s,x,\xi)$ in the  standard Frechet topology on $\Ci(\Cf\times
U,(\Rf^N)^{*}\otimes\Rf^N)$. So we obtain
\begin{eqnarray*}
\left.\int_{|\xi|\geq 1} \as_n\si(x,\xi) \ \hat{d}\xi\right|^{{\rm
mer}} & = & \frac{1}{\a s}\int_{|\xi| = 1} \left( \ \Is_n(x,\xi)
- s\,(\log\as)_n(x,\xi) + o(s) \ \right) \ \hat{d}\xi \\
& = & -\frac{1}{\a}\int_{|\xi| = 1} (\log\as)_n(x,\xi) \,
\hat{d}\xi + o(s^0) .
\end{eqnarray*}
Hence
\begin{equation}\label{e:jequalton}
\left.\int_{|\xi|\geq 1} \as_n\si(x,\xi) \
\hat{d}\xi\right|_{s=0}^{{\rm mer}} = -\frac{1}{\a}\int_{|\xi| =
1} (\log\as)_n(x,\xi) \ \hat{d}\xi \ .
\end{equation}

\vskip 2mm

For general $s\in\Cf$ the first (non-homogeneous) term on the
right side of \eqref{e:twoparts} is a more complicated expression
leading to global poles. However, at $s=0$ this is local and, from
\eqref{e:symbolpowerzero}, given by
\begin{equation*}
    \left.\int_{|\xi|\leq 1} \as_j\si(x,\xi) \
\hat{d}\xi\right|^{{\rm mer}}_{s=0} = \int_{|\xi|\leq 1}
\Is_j(x,\xi) \ \hat{d}\xi \ .
\end{equation*}

\noi Hence
\begin{eqnarray}
\left.\int_{|\xi|\leq 1} \as_j\si(x,\xi) \
\hat{d}\xi\right|_{s=0}^{{\rm
mer}} & = & 0  \hskip 10mm j\neq 0 \ . \label{e:firsttermjnotzero} \\
\left.\int_{|\xi|\leq 1} \as_0\si(x,\xi) \
\hat{d}\xi\right|_{s=0}^{{\rm mer}} & =& \frac{1}{(2\pi)^n} \, .
\, {\rm vol}(B^n)   \ , \label{e:firsttermjzero}
\end{eqnarray}
where $B^n$ is the $n$-ball.

\vskip 2mm

Thus from \eqref{e:slarge}, \eqref{e:globalbit},
\eqref{e:twoparts}, \eqref{e:jnotzeronorn}, \eqref{e:jzero},
\eqref{e:jequalton}, \eqref{e:firsttermjnotzero},
\eqref{e:firsttermjzero} we have
\begin{eqnarray*}
\left. K(A_{\th}\si,x,x)\right|_{s=0}^{{\rm mer}} & = &
-\int_{\Rf^n} \s(\Pi_0(A))(x,\xi) \ \hat{d}\xi dx \\
& & -\frac{1}{\a}\int_{|\xi| =1} (\log\as)_n(x,\xi) \,
\hat{d}S(\xi) - \, \frac{1}{(2\pi)^n}.\frac{{\rm vol}(S^{n-1})}{n}
+ \frac{1}{(2\pi)^n} \, . \, {\rm vol}(B^n)  \\ & = &
-\int_{\Rf^n} \s(\Pi_0(A))(x,\xi) \ \hat{d}\xi dx -
\frac{1}{\a}\int_{|\xi|=1} (\log\as)_n(x,\xi) \, \hat{d}S(\xi) \ .
\end{eqnarray*}

\noi Hence

$$\int_M\int_{|\xi|=1} \tr((\log\as)_n(x,\xi)) \, \hat{d}S(\xi) \, dx
=$$
\begin{equation*}
-\a\left(\int_M \left.\tr(K(A_{\th}\si,x,x)\right|_{s=0}^{{\rm
mer}}) \, dx  + \int_M\int_{\Rf^n} \tr(\s(\Pi_0(A))(x,\xi)) \
\hat{d}\xi dx\right)   \ ,
\end{equation*}

\noi that is,

\begin{equation*}
\res(\log\,A) = -\a\,\left(\,\z(A,0)|^{{\rm mer}} +
\Tr(\Pi_0(A))\,\right) \ .
\end{equation*}
\end{proof}

%\begin{eqnarray}
%\z(A,s)|^{{\rm mer}}_{s=0} & := & \int_M \int_{\Rf^n}
%\tr(K(A\si_{\th},x,x)|^{{\rm mer}}_{s=0}) \ \hat{d}\xi dx \nonumber \\
%& = & -\int_M\int_{\Rf^n} \tr(\s(\Pi_0(A))(x,\xi)) \ \hat{d}\xi dx
%  \ + \ \sum_{j=0}^{n}
%\int_M\int_{\Rf^n} \tr(\as_{\th,j}\si(x,\xi)|^{{\rm mer}}_{s=0}) \ \hat{d}\xi dx  \nonumber \\
%& = & - \, \Tr(\Pi_0(A))   \ + \ \sum_{j=0}^{n} \int_M\int_{\Rf^n}
%\tr(\as_{\th,j}\si(x,\xi)|^{{\rm mer}}_{s=0}) \ \hat{d}\xi dx
%\label{e:localbits}
%\end{eqnarray}

\vskip 5mm

\noi {\bf Proof of \thmref{t:idpluslowerorder}}

\vskip 5mm

\begin{proof}
We have,
\begin{equation*}
\log(I + \Qs)  = \frac{i}{2\pi}\int_{\Gamma_{\theta}}\log\la \ (I
+ \Qs -\la I)\ii \ d\la \ ,
\end{equation*}
where the finite contour $\Gamma_{\theta}$ encloses, in particular
$1$. Iterating
\begin{equation*} (I + \Qs -\la I)\ii =(1-\la)\ii I\ -(1-\la)\ii
\Qs(I + \Qs -\la I)\ii
\end{equation*}
yields
\begin{equation*}
(I + \Qs -\la I)\ii = \sum_{j=0}^m (-1)^j(1-\la)^{-j-1}\Qs^{j} +
(-1)^{m+1}(1-\la)^{m} \Qs^m(I + \Qs -\la I)\ii \ ,
\end{equation*}
and so
\begin{equation}\label{e:resexp}
\log(I + \Qs)  = \sum_{j=0}^m (-1)^j \, \Qs^{j} \,
\frac{i}{2\pi}\int_{\Gamma_{\theta}} \log\la \,(1-\la)^{-j-1} \
d\la + R(\Qs,m) \ ,
\end{equation}
where
\begin{equation*}
R(\Qs,m)  = (-1)^{m+1} \Qs^m \
\frac{i}{2\pi}\int_{\Gamma_{\theta}}\log\la \ (1-\la)^{m} (I + \Qs
-\la I)\ii \ d\la \ .
\end{equation*}
$R(\Qs,m)$ is a classical $\pdo$ of order $-km$ and so for any
positive integer $m$ with $m k < -n$ we have $\res(R(\Qs,m))=0$.
All operators in \eqref{e:resexp} are integer order and so we can
use the linearity of $\res$ in \lemref{lem:linearity} to find
\begin{equation*}
\res(\log(I + \Qs))  = \sum_{j=1}^{\left[\frac{n}{|k|}\right]}
(-1)^j \, \res(\Qs^{j}) \,
 \frac{i}{2\pi}\int_{\Gamma_{\theta}}\log\la
\,(1-\la)^{-j-1} \ d\la  \ ,
\end{equation*}
the summation beginning now from $j=1$, since $\res(I) =0$ and the
contour integral is zero for $j=0$. Since $\Gamma_{\theta}$
encloses $1$, then for $j\geq 1$
\begin{equation*}
\frac{i}{2\pi} \int_{\Gamma_{\theta}}\log\la \ (1-\la)^{-j-1} \
d\la  = \frac{1}{j}\ \frac{i}{2\pi} \int_{\Gamma_{\theta}}\la\ii \
(1-\la)^{-j} \ d\la = -\frac{1}{j}
\end{equation*}
and  we reach the conclusion.
\end{proof}

\vskip 5mm

\noi {\bf Proof of \thmref{t:zetapolynomial}}

\vskip 5mm

\begin{proof}
Let $A \in\Psi^{\a}(E)$ be an elliptic $\pdo$, admitting a
principal angle. Let $Q$ be a parametrix for $A$, so that
\begin{equation}\label{e:2sidedparametrix}
AQ-I = s_{\o}\in \Psi^{-\infty}(E) \ , \hskip 15mm QA-I =
\tilde{s}_{\o}\in \Psi^{-\infty}(E)
\end{equation}
are smoothing operators. For any smoothing operator $\k_{\o}\in
\Psi^{-\infty}(E)$ one has by \corref{cor:detresstable}
\begin{equation}\label{e:stable}
{\rm det}_{{\rm res}}(A + \k_{\o}) = {\rm det}_{{\rm res}}(A) \ .
\end{equation}

\noi Let $B \in\Psi^{\a}(E)$ with $\a-\b\in\Nf$. Then by
\eqref{e:detresmult}, which from the proof of
\thmref{e:detresmult} is seen to hold logarithmically,
\begin{eqnarray*}
\log{\rm det}_{\res} (A + B) & = &\log{\rm det}_{\res}
(AQ-s_{\o})(A + B) \\
& = &\log{\rm det}_{\res} (A Q A + A Q B + t_{\o}) \\
& = &\log{\rm det}_{\res} (A Q A + A Q B) \\
& = & \log{\rm det}_{\res}A  +\log{\rm det}_{\res} (Q A + Q
B)\\
& = & \log{\rm det}_{\res}A  +\log{\rm det}_{\res} (I + Q B) \ ,\\
\end{eqnarray*}
where $t_{\o}\in\Psi^{-\infty}(E)$ and for the final equality we
use \eqref{e:2sidedparametrix} and \eqref{e:stable}. Rewriting in
terms of \eqref{e:resdetzetazero1} and \eqref{e:degreezero} this
reads
\begin{equation*}
-\a \,(\,\z(A+B,0)|^{{\rm mer}} + \h_0(A+B)\,)= -\a
\,(\,\z(A,0)|^{{\rm mer}}  + \h_0(A)\,)+ \sum_{j=0}^M
\frac{(-1)^j}{j} \, \res\,\left( \,(QB)^j \, \right) \ .
\end{equation*}
The sum terminates when $\ord(QB).j < -n$, so we may take
\begin{equation*}
M = \left[\frac{n}{\a-\b} \right] \ .
\end{equation*}
Replacing $B$ by $B[t] = B_0 + B_1 \, t + \ldots B_d \, t^d$ now
proves \eqref{e:zetapolynomial1}.
\end{proof}

\vskip 8mm

\noi {\small King's College, \newline \noi London. \vskip 2mm \noi
Email: sgs@mth.kcl.ac.uk}


\begin{thebibliography}{99}

\bibitem[Bo]{Bo} J. Bost, {\em Fibr$\acute{{\rm e}}$s d$\acute{{\rm e}}$terminants,
d$\acute{{\rm e}}$terminants r$\acute{{\rm e}}$gularises et
mesures sur le espaces de modules de courbes complexes}, Asterique
{\bf 152} (1988), 113--149.

\bibitem[BrOr]{BrOr} T. Branson and B. \O rsted, {\em Conformal geometry and global
invariants}, Diff. Geom. Appl. {\bf 1} (1991), 279--308.

\bibitem[Bu]{Bu} T. Burak, {\em On spectral projections of elliptic operators},
Ann. Scuola Norm. Sup. Pisa {\bf 24} (1970), 209--230.

\bibitem[Gi]{Gi} P. Gilkey, {\em Invariance Theory, the heat
equation and the Atiyah-Singer Index Theorem}.  2nd Edition, CRC
Press, 1995.

\bibitem[Gr]{Gr} G. Grubb, {\em A resolvent approach to traces
and zeta Laurent expansions}, AMS Contemp. Math. Proc., vol. 366,
2005, pp.67--93, arXiv: math.AP/0311081.

\bibitem[Gr2]{Gr2}  G. Grubb, {\em On the logarithmic component
in trace defect formulas}, preprint 2004, arXiv: math.AP/0411483.

\bibitem[GrSe]{GS} G. Grubb and  R. Seeley, {\em Weakly parametric
pseudodifferential operators and Atiyah-Patodi-Singer boundary
problems},  Invent. Math. {\bf 121} (1995), 481-529.

\bibitem[Gu]{Gu} V. Guillemin,
{\em A new proof of Weyl's formula on the asymptotic distribution
of eigenvalues}, Adv. Math. {\bf 102} (1985), 184--201.

\bibitem[Ja]{J} N. Jacobson, {\em Lie Algebras}, Interscience
Tracts in Pure and Appl. Math. {\bf 10} 1962. Wiley, New York.

\bibitem[Ka]{Ka} C. Kassel, {\em Le r$\acute{{\rm e}}$sidu
non-commutatif [d'Apr$\grave{{\rm e}}$s M.Wodzicki]},  Asterique
{\bf 177-178} (1989), 199--229.

\bibitem[KoVi]{KV} M. Kontsevich and S. Vishik,
{\em Determinants of elliptic pseudodifferential operators} arXiv:
hep-th/9404046 (1994); {\em Geometry of determinants of elliptic
operators}, Funct. Anal. on the Eve of the 21st Century {\bf 1},
Birkhauser, Progr. Math. {\bf 131}, 1995, pp.173--197.

\bibitem[McSi]{McSi} H. McKean and I. Singer,
{\em Curvature and the eigenvalues of the Laplacian}, J. Diff.
Geom. {\bf 1} (1967), 43--69.

\bibitem[Mu]{Mu} W. M$\ddot {\rm u}$ller,
{\em Relative zeta functions, relative determinants and scattering
theory}, Comm. Math. Phys. {\bf 192} (1998), 309--347.

\bibitem[Ok1]{Ok1} K. Okikiolu, {\em The Campbell-Hausdorff theorem
for elliptic operators and a related trace formula}, Duke Math. J.
{\bf 79} (1995), 687--722.

\bibitem[Ok2]{Ok2} K. Okikiolu, {\em The multiplicative anomaly for
determinants of elliptic operators}, Duke Math. J. {\bf 79}
(1995), 723--750.

\bibitem[Pa]{Pa} S. Paycha, Anomalies and regularization techniques in mathematics and
physics, preprint (Colombia) 2004.

\bibitem[PaRo]{PaRo} S. Paycha and S. Rosenberg,
{\em Traces and characteristic classes on loop groups}, in
Infinite dimensional groups and manifolds,  ed. V.Turavev,T.
Wurzbacher, de Gruyter, Berlin, 2002.

\bibitem[PaSc]{PaSc} S. Paycha and S. Scott, {\em The Laurent expansion
for regularized integrals of holomorphic symbols}, preprint 2004.

\bibitem[Po]{Po} R. Ponge, {\em Spectral asymmetry, zeta functions, and the noncommutative residue},
Preprint, arXiv: math.DG/0310102 (2005).

\bibitem[Sc]{Sc} S. Scott, {\em Zeta determinants on manifolds with
boundary}, J. Funct. Anal. {\bf 192} (2002), 112--185.

\bibitem[ScZa]{ScZa} S. Scott and D. Zagier, {\em A symbol proof of the
local index theorem}, preprint 2004.

\bibitem[Se]{Se} R. T. Seeley, {\em Complex powers of an elliptic
operator}, AMS Proc. Symp. Pure Math. X, 1966, AMS Providence,
1967, pp. 288--307.

\bibitem[Sh]{Sh} M. A. Shubin, {\em Pseudodifferential Operators
and Spectral Theory}, Springer, 1987.

\bibitem[Si]{Si} B. Simon, {\em Trace Ideals and their Applications},
LMS Lecture Notes {\bf 35}, CUP, 1979.

\bibitem[Wo1]{Wo1} M. Wodzicki, {\em Spectral asymmetry and zeta functions},
Invent. Math. {\bf 66} (1982), 115--135.

\bibitem[Wo2]{Wo2} M. Wodzicki, {\em Local invariants of spectral
asymmetry}, Invent. Math. {\bf 75} (1984), 143--178.

\bibitem[Wo3]{Wo3} M. Wodzicki, {\em Non-commutative residue,
Chapter I. Fundamentals}, K-Theory, Arithmetic and Geometry,
Springer Lecture Notes {\bf 1289}, 1987, pp.320--399.

\end{thebibliography}
\end{document}